\pdfoutput=1
\documentclass[final]{siamart1116}
\usepackage{url}
\usepackage{amssymb}
\usepackage{amsmath}
\usepackage{multirow}
\usepackage{booktabs}
\usepackage[numbers,square,sort&compress]{natbib}
\usepackage{graphicx}
\usepackage{subcaption}

\usepackage{pythonhighlight}

\usepackage{paralist}

\DeclareMathOperator{\Div}{div}
\DeclareMathOperator{\Grad}{grad}
\newcommand{\inner}[2]{\langle #1, #2 \rangle}
\newcommand{\R}{\mathbb{R}}

\newcommand{\ddt}[1]{{#1}_t}
\newcommand{\dx}{\, \mathrm{d} x}
\newcommand{\ds}{\, \mathrm{d} s}

\newcommand{\dP}{\, \mathrm{d} P}

\newcommand{\triang}{\mathcal{T}}

\newcommand{\timestep}{\kappa}
\newcommand{\version}{2016.1}
\newcommand{\avg}[1]{\langle #1 \rangle}

\newtheorem{example}{Example}

\title{{Automated adjoints of coupled PDE-ODE systems}
  \thanks{\today{}.
    \funding{This research is supported by a Center of Excellence
      grant awarded to the Center for Biomedical Computing at Simula
      Research Laboratory from the Research Council of Norway, by
      EPSRC grants EP/K030930/1 and EP/M011151/1, a NOTUR grant
      NN9316K and the generous support of Sir Michael Moritz and
      Harriet Heyman.}}}
\date{\today}
\author{
  P. E. Farrell\thanks{Mathematical Institute, University of Oxford, OX2 6GG, UK (\email{patrick.farrell@maths.ox.ac.uk})}
  \and
  J. E. Hake\thanks{Ski Videregående skole, Ski, Norway} (\email{johan.hake@gmail.com})
  \and
  S. W. Funke\thanks{Simula Research Laboratory, P.O. Box 134, 1325 Lysaker, Norway} (\email{simon@simula.no})
  \and
  M. E. Rognes\thanks{Simula Research Laboratory, P.O. Box 134, 1325 Lysaker, Norway} (\email{meg@simula.no})
}

\begin{document}

\maketitle

\begin{abstract}
Mathematical models that couple partial differential equations (PDEs)
and spatially distributed ordinary differential equations (ODEs) arise
in biology, medicine, chemistry and many other fields. In this paper
we discuss an extension to the FEniCS finite element software for
expressing and efficiently solving such coupled systems. Given an ODE
described using an augmentation of the Unified Form Language (UFL) and
a discretisation described by an arbitrary Butcher tableau, efficient
code is automatically generated for the parallel solution of the
ODE. The high-level description of the solution algorithm also
facilitates the automatic derivation of the adjoint and tangent
linearization of coupled PDE-ODE solvers. We demonstrate the
capabilities of the approach on examples from cardiac
electrophysiology and mitochondrial swelling.
\end{abstract}



\section{Introduction}
\label{sec:introduction}

In this paper we discuss solvers for systems involving spatially
dependent ordinary differential equations (ODEs) of the following
form: given an initial condition $y_0(x)$, find $y = y(x, t)$ such
that
\begin{align}
  \label{eq:intro:ivp}
  \ddt{y}(x, t) = f(y, x, t), \quad y(x, t_0) = y_0(x),
\end{align}
for all $x$ in a point set $X \subset \Omega \subseteq \R^d$, where
$d$ is the spatial dimension. The subscript $t$ refers to
differentiation in time.  The right-hand side $f$ cannot depend on
spatial derivatives of $y$; the ODE is decoupled at different points.
Problems of this form often arise when discretizing time-dependent
mathematical models that couple PDEs with spatially distributed
systems of ODEs via operator splitting. Examples of application areas
include cardiac electrophysiology in
general~\citep{SundnesEtAl2006,Niederer2011} and cardiac ion channel
modeling in particular~\citep{FitzHugh1961}, mitochondrial
swelling~\citep{Eisenhofer2013}, groundwater flow and
contamination~\citep{WangEtAl2003}, pulmonary gas
transport~\citep{WhiteleyEtAl2002, BignellJohnston2007}, and
plasma-enhanced chemical vapor deposition~\citep{Geiser2011}.
We also discuss the automated derivation of the adjoint and tangent
linearization of such models:
these
can be used to identify the sensitivity of the solution to model parameters,
solve inverse problems for unknown parameters, and characterize
the stability of trajectories.

For concreteness, we present two biological examples of coupled
PDE-ODE systems where problems of the form \eqref{eq:intro:ivp}
arise. We return to numerical results for these examples in
Section~\ref{sec:applications}.

\begin{example}[The bidomain equations]
As our first motivating example, we will consider the bidomain
equations for the propagation of an electrical signal in a
non-deforming domain $\Omega$~\citep{SundnesEtAl2006}: find the
transmembrane potential $v = v(x, t)$, the extracellular potential
$u_e(x, t)$ and additional state variables $s = s(x, t)$ such that for
$t \in (0, T]$:
\begin{subequations}
  \label{eq:bidomain}
\begin{align}
  \label{eq:bidomain:1}
  \ddt{v} - \Div (M_i \Grad v + M_i \Grad u_e) &= - I_{\rm ion}(v, s) + I_s
  &\text{in } \Omega, \\
  \label{eq:bidomain:2}
  - \Div \left (M_i \Grad v + (M_i + M_e) \Grad u_e \right ) &= 0
  &\text{in } \Omega, \\
  \label{eq:bidomain:3}
  \ddt{s} &= F(v, s)
  &\forall x \in \Omega.
\end{align}
\end{subequations}
In~\eqref{eq:bidomain}, $I_{\rm ion}$ is a given nonlinear function
describing ionic currents and $F$ defines a system of nonlinear
functions, while $M_i$ and $M_e$ are the intracellular and
extracellular conductivity tensors, respectively, and $I_s$ is a given
stimulus current. The function $F$ cannot depend on spatial
derivatives of $v$ and $s$; it defines a pointwise system of ODEs. The
specific form of $I_{\rm ion}$ and $F$ are typically prescribed by a
given \emph{cardiac cell model} and may vary greatly in complexity:
from involving a single state variable $s$ such as the FitzHugh-Nagumo
model~\citep{FitzHugh1961} to models with e.g.~$41$ state variables
such as the model of O'Hara et al.~\citep{OHara2011}. The
system~\eqref{eq:bidomain} is closed with appropriate initial and
boundary conditions. After the application of operator splitting, the
ODE step is decoupled from the PDE step and is a system of the form of
\eqref{eq:intro:ivp}.
\end{example}

\begin{example}[Mitochondrial swelling]
As a second example, we consider a model proposed
in~\citep{Eisenhofer2013} to describe the swelling of
mitochondria. Mitochondrial swelling plays a key role in the process
of programmed cell death (\emph{apoptosis}) and thus for the life
cycle of cells. Mathematically, we consider the following
model~\citep[p. 26]{Eisenhofer2013}: find the calcium concentration $u
= u(x, t)$ and the densities of mitochondria states $N_i = N_i(x, t)$
for $i = 1, 2, 3$ such that
\begin{subequations}
  \label{eq:mitochondria}
  \begin{align}
    \label{eq:mito:1}
    \ddt{u} &= d_1 \Delta (|u|^{q-2} u) + d_2 g (u) N_2
    &\text{in } \Omega, \\
    \label{eq:mito:2}
    \ddt{N_1} &= - f(u) N_1
    &\forall x \in \Omega, \\
    \label{eq:mito:3}
    \ddt{N_2} &= f(u) N_1 - g(u) N_2
    &\forall x \in \Omega, \\
    \label{eq:mito:4}
    \ddt{N_3} &= - g(u) N_2
    &\forall x \in \Omega.
  \end{align}
\end{subequations}
where $d_1 \geq $ is a diffusion coefficient, $d_2 \geq 0$ is a
feedback parameter, $q \ge 2 \in \mathbb{N}$ determines the
nonlinearity of the diffusion-type operator, and $g$ and $f$ are
prescribed functions. The functions $N_1$, $N_2$ and $N_3$
describe the densities of unswollen, swelling and completely
swollen mitochondria respectively. Equation \eqref{eq:mito:1} is a
spatially-coupled PDE, while equations
\eqref{eq:mito:2}-\eqref{eq:mito:4} define pointwise ODEs. The system
is closed with initial conditions and Dirichlet boundary
conditions for $u$. Again, after operator splitting, a subproblem of
the form \eqref{eq:intro:ivp} results.
\end{example}

Over the last decade, there has been a growing interest in
computational frameworks for the rapid development of numerical solvers
for PDEs such
as the FEniCS Project~\citep{LoggMardalEtAl2012a}, the Firedrake
Project~\citep{Rathgeber2015}, and
Feel++~\citep{PrudhommeEtAl2012}. The rapid development of
solvers is achieved by offering high-level abstractions
for the expression of such problems. Each of these projects
provides a domain-specific language for specifying finite element
variational formulations of PDEs and associated software for their
efficient solution. The methods presented in this paper are
implemented in the the FEniCS Project and
dolfin-adjoint~\citep{FarrellEtAl2013}; dolfin-adjoint automatically
derives discrete adjoint and tangent linear models from a FEniCS
forward model. These allow for the efficient computation of functional
gradients and Hessian-vector products, and are essential ingredients in stability
analyses, parameter identification and inverse problems.

However, these systems do not currently efficiently extend to the kind
of coupled PDE-ODE systems arising in computational biology.  While
monolithic finite element discretizations of coupled PDE-ODE systems
such as~\eqref{eq:bidomain} or~\eqref{eq:mitochondria} may easily be
specified and solved using the current software features in FEniCS,
this approach will involve an unreasonably large computational
expense, as the monolithic system is highly nonlinear and the solver
cannot exploit the fact that the ODE is spatially decoupled.

Operator splitting is the method of choice for such coupled PDE-ODE
systems~\citep{SundnesEtAl2006}, and is implemented as standard in
many hand-written codes, see
e.g.~\citep{VigmondEtAl2008,MiramsEtAl2013} in the context of cardiac
electrophysiology. Operator splitting decouples a PDE-ODE system into
a spatially-coupled system of PDEs and a spatially-decoupled
collection of ODE systems. This collection of ODE systems then
typically takes the form~\eqref{eq:intro:ivp}, which can be solved
using well-known temporal discretization methods. Heretofore it has
not been possible to specify spatially-decoupled ODE systems in
high-level PDE frameworks such as the FEniCS or Firedrake projects.

This work addresses the gap in available abstractions, algorithms
and software for arbitrary PDE-ODE systems.
We introduce high-level
domain-specific language constructs for specifying collections of ODE
systems and for specifying multistage ODE schemes via Butcher
tableaux. This has several major benefits. By automatically generating
the solver from a high-level description of the problem, practitioners
can flexibly explore a range of models; this is especially important
in biological problems where the model itself is uncertain.  Another
advantage is that a high-level description facilitates the
automated derivation of the associated tangent linear and adjoint
models.  This represents a significant saving in the time taken to
investigate questions of biological interest.

The main new contributions of this paper are (i) the extension of the
FEniCS finite element system to enable efficient the large-scale
forward solution of coupled, time-dependent PDE-ODE systems via
operator splitting, and (ii) the extension of dolfin-adjoint to
automatically derive and solve the associated tangent linear and
adjoint models, enabling efficient automated computation of functional
gradients for use in e.g.~optimization or adjoint-based sensitivity
analysis.

This paper is organized as follows. In
Section~\ref{sec:problemsetting}, we describe the general operator
splitting setting, the resulting separate PDE and ODE systems, and
natural discretizations of these. We continue in
Section~\ref{sec:ode_schemes} by briefly describing the
well-established ODE schemes that we consider in this work: the
multistage and Rush-Larsen families. We discuss the adjoint and
tangent linear discretizations of a general operator splitting scheme
and derive the adjoint and tangent linear models for the multistage
and Rush-Larsen schemes in Section~\ref{sec:odes:adjoints}. Key
features of our new implementation are described in
Section~\ref{sec:implementation}. We present numerical results for two
different application examples in Section~\ref{sec:applications}, and
use these to evaluate the performance of our implementation. In
Section~\ref{sec:conclusion}, we provide some concluding remarks and
discuss current limitations and possible future extensions.

\section{Operator splitting for coupled PDE-ODE systems}
\label{sec:problemsetting}

\subsection{Operator splitting}
\label{sec:operator_splitting}

A classical approach to the solution of the bidomain
system~\eqref{eq:bidomain}~\citep{SundnesEtAl2006} and similar systems
is to apply operator splitting. This approximately solves the full
system of equations by alternating between the solution of a system of
PDEs and a system of ODEs defined over each time interval. The
advantage of this approach is that it decouples the solution of the
typically highly nonlinear ODEs from the spatially coupled PDEs at
each time step. We will consider this general setting, but
use~\eqref{eq:bidomain} as a concrete example.

Applying a variable order operator split to~\eqref{eq:bidomain}, the
resulting scheme reads: given initial conditions $v^0$, $s^0$, an order
parameter $\theta \in [0, 1]$, and time points $\{t_0, \dots, t_{N}\}$
with an associated timestep $\timestep_n = t_{n+1} - t_{n}$, then for
each time step $n = 0, 1, \dots, N - 1$:
\begin{compactenum}
\item
  Compute $v^\ast$ and $s^\ast$ by solving the ODE system
  \begin{subequations}
    \label{eq:bidomain:odes}
    \begin{align}
      \ddt{v} &= - I_{\rm ion} (v, s) \\
      \ddt{s} &= F(v, s)
    \end{align}
  \end{subequations}
  over $\Omega \times [t_n, t_n + \theta \kappa_n]$ with
  initial conditions $v^{n}, s^{n}$.
\item
  Compute $v^{\dagger}$ and $u^{n+1}_e$ by solving the PDE system
  \begin{subequations}
    \label{eq:bidomain:pdes}
    \begin{align}
      \ddt{v} - \Div (M_i \Grad v + M_i \Grad u_e) &= I_s \\
      - \Div (M_i \Grad v + (M_i + M_e) \Grad u_e) &= 0
    \end{align}
  \end{subequations}
  over $\Omega \times [t_n, t_{n+1}]$ with initial condition
  $v^{\ast}$.
\item
  If $\theta < 1$, compute solutions $v^{n+1}$ and $s^{n+1}$
  solving~\eqref{eq:bidomain:odes} over $\Omega \times [t_n + \theta
    \kappa_n, t_{n+1}]$ with initial conditions $v^{\dagger}$ and
  $s^{\ast}$.
\end{compactenum}
The split scheme relies on the repeated solution of a
nonlinear system of ODEs and (in this case) a linear system of
PDEs. The main advantage of the approach is that it allows for the
separate discretization and solution of the ODEs and PDEs, with the
respective solutions as feedback into the other system.
For $\theta = 1/2$, the resulting Strang splitting scheme is
second-order accurate; for other values of $\theta$ the
resulting scheme is first-order accurate.

\subsection{Discretization of the separate PDE and ODE systems}

Suppose now that the relevant PDE system
(e.g.~\eqref{eq:bidomain:pdes}) is discretized in space by a finite
element method defined over a mesh $\triang_h$ of the domain $\Omega$,
and in time by some suitable temporal discretization. Efficient
solution algorithms for such discretizations are well-established
(such as~\citep{LoggMardalEtAl2012a}) and will not be detailed further
here.

At each iteration the solution of the PDE system relies on the
solution of the ODE system (e.g.~\eqref{eq:bidomain:odes}), or at
least the ODE solution evaluated at some finite set of points $X =
\{x_i\}_{i=1}^{|X|}$ in space. For instance, the set of points $X$ may
be taken as the nodal locations of the finite element degrees of
freedom, or the quadrature points of the mesh. As a consequence of this and the
spatial locality of ODEs, a natural approach to discretizing the
system of ODEs is to step the ODE system forward in time at this set
of points $X$. Typically, $|X|$ is very large and the efficient
repeated solution of these systems of ODEs is key. This is the setting
that we focus on next.

\section{Solution schemes for the ODE systems}
\label{sec:ode_schemes}

The initial value problem~\eqref{eq:intro:ivp} decouples in space, and
henceforth we consider its solution for a fixed $x$. With a minor
abuse of notation, let $y = y(x)$ and $f(y, t) = f(y, x, t)$ so
that~\eqref{eq:intro:ivp} reads as the classical ODE problem: find $y
\in C^1([T_0, T_1]; \mathbb{R}^m)$ for a certain $m \in \mathbb{N}$
such that
\begin{align}
  \label{eq:ivp}
  \ddt{y}(t) = f(y, t), \quad y(T_0) = y_0.
\end{align}
for $t \in [T_0, T_1]$. There exists a wide variety of solution
schemes for~\eqref{eq:ivp} including but not limited to multistage,
multistep, and IMEX schemes, see e.g.~\citep{Butcher2008}. In this
work, we focus on two classes of schemes that are widely used in
computational biology: multistage schemes and so-called Rush-Larsen
schemes. The Rush-Larsen schemes are commonly used in cardiac
electrophysiology in general and for discretizations of the bidomain
equations~\eqref{eq:bidomain} in particular. These two classes of
schemes are detailed below. Keeping the general iterative operator
splitting setting in mind, we present the schemes on a single
time-step $[T_0, T_1]$ with $\kappa = T_1 - T_0$ for brevity of
notation.

\subsection{Multistage schemes}
\label{sec:multistage:schemes}

An $s$-stage multistage scheme for~\eqref{eq:ivp} is defined by a set
of coefficients $a_{ij}$, $b_i$ and $c_i$ for $i, j = 1, \dots, s$,
commonly listed in a so-called \emph{Butcher tableau}; for more
details, see e.g.~\citep{Butcher2008}. Given $y^0$ at $T_0$, the
scheme finds the stage variables $k_i$ for $i = 1, \dots, s$
satisfying
\begin{equation}
  \label{eq:multistage:stage}
  k_i = f(y^{0} + \kappa \sum_{j=1}^s a_{ij} k_j, T_0 + c_i \kappa),
\end{equation}
and subsequently sets the solution $y^1$ at $T_1$ via
\begin{equation}
  \label{eq:multistage:final}
  y^{1} = y^{0} + \kappa \sum_{i=1}^s b_i k_i.
\end{equation}
Note that~\eqref{eq:multistage:stage} defines a system of (non-linear)
equations to solve for the stage variables $k_i$ ($i = 1, \dots, s$)
if $a_{ij} \not = 0$ for any $j \geq i$. Our implementation demands
that $a_{ij} = 0$ for $j > i$, i.e.~does not allow the computation of
earlier stages to depend on the values of later stages.

\subsection{The Rush-Larsen scheme and its generalization}
Recall that $y = \{y_i\}_{i=1}^m$ and $f(y) = \{f_i(y)\}_{i=1}^m$. The
original Rush-Larsen scheme employs an exponential integration scheme
for all linear terms for each component $f_i$ (for $i = 1, 2, \dots,
M$), and a forward Euler step for all non-linear
terms~\citep{RushLarsen1978}. Let $J_{i} = \frac{\partial
  f_i}{\partial y_i}(y^0, T_0)$ be the i'th diagonal component of the
Jacobian at time $T_0$ for $i = 1, \dots, M$, and assume that $J_{i}$
is non-zero. For brevity, denote $f_i(y^0, T_0) = f_{i}^0$. The
Rush-Larsen scheme then computes $y_i^{1}$ at $T_1$ by
\begin{equation}
  \label{eq:rl1}
  y_i^{1} = y_i^{0} + \left\{
    \begin{array}{ll}
      J_{i}^{-1} f_{i}^0 \left(e^{\kappa J_{i}}-1 \right) & \mbox{if $f_i$ is linear in $y_i$} , \\
      \kappa f_{i}^0 & \mbox{if $f_i$ is not linear in $y_i$}.
    \end{array}
\right .
\end{equation}%

For stiff systems of ODEs, the forward Euler step in~\eqref{eq:rl1}
can become unstable for large $\kappa$. This motivates a generalized
version of the Rush-Larsen scheme~\citep{SundnesEtAl2009}. In this
generalization, the exponential integration step is used for all
components of $y$ reducing the scheme to
\begin{equation}
  \label{eq:grl1}
  y_i^{1} = y_i^{0} + J_{i}^{-1} f_i^0 \left(e^{\kappa J_{i} }-1 \right)
\end{equation}
Both~\eqref{eq:rl1} and~\eqref{eq:grl1} are first order accurate in
time~\citep{RushLarsen1978, SundnesEtAl2009}.

Both the original and the generalized Rush-Larsen schemes can be
developed into second order schemes by repeated use as follows. The
solution $y^{\frac{1}{2}}$ at $T_{\frac{1}{2}} = T_0 +
\frac{\kappa}{2}$ is computed in the first step and is used in $f$ and
its linearization $J$ to compute $y^{1}$ in the second step. Let
$f_{i}^{\frac{1}{2}} = f_i(y^{\frac{1}{2}}, T_{\frac{1}{2}})$, and let
$J_{i, \frac{1}{2}} = \frac{\partial f_i}{\partial
  y_i}(y^{\frac{1}{2}}, T_{\frac{1}{2}})$ be the diagonal component of
the Jacobian at time $T_{\frac{1}{2}}$ for $i = 1, \dots, M$. More
precisely, the second order version of~\eqref{eq:rl1}
and~\eqref{eq:grl1} is then given by
\begin{subequations}
  \label{eq:rl2}
  \begin{align}
    y_i^{\frac{1}{2}} &= y_i^{0} + \left\{
    \begin{array}{ll}
      J_{i}^{-1} f_i^0 \left(e^{\frac{\kappa}{2}  J_{i}} - 1 \right) & \mbox{if $f_i$ is linear in $y_i$}\\
      \frac{\kappa}{2} f_i^0 & \mbox{if $f_i$ is not linear in $y_i$}
    \end{array}
    \right.\\
    y_i^{1} &= y_i^{0} + \left\{
    \begin{array}{ll}
      J_{i, \frac{1}{2}}^{-1} f_i^{\frac{1}{2}} \left (e^{\kappa J_{i, \frac{1}{2}}} - 1 \right) & \mbox{if $f_i$ is linear in $y_i$}\\
      \kappa f_{i}^{\frac{1}{2}} & \mbox{if $f_i$ is not linear in $y_i$}
    \end{array}
    \right.
  \end{align}%
\end{subequations}%
for the original Rush-Larsen scheme and
\begin{subequations}
  \label{eq:grl2}
  \begin{align}
    y_i^{\frac{1}{2}} &= y_i^{0} + J_{i}^{-1} f_i^0  \left(e^{\frac{\kappa}{2} J_{i}}-1 \right) \\
    \label{eq:grl2_1}
    y_i^{1} &= y_i^{0} + J_{i, \frac{1}{2}}^{-1} f_{i}^{\frac{1}{2}} \left( e^{\kappa J_{i, \frac{1}{2}}} - 1 \right)
  \end{align}
\end{subequations}
for the generalized Rush-Larsen scheme. Here, \eqref{eq:grl2_1} is
simplified compared to the scheme in~\citep{SundnesEtAl2009}: instead
of evaluating $f$ and $J$ at $\bar{y}^{(i)}$ (cf. eq.~(7)
in~\citep{SundnesEtAl2009}) we use $y_{i}^{\frac{1}{2}}$. 

\section{Adjoints and tangent linearizations of coupled PDE-ODE systems}
\label{sec:odes:adjoints}
Our aim is to automatically derive the adjoint and tangent linear
equations of operator splitting schemes in a manner that allows for
efficient solution of the resulting systems. We proceed as follows: we
first state the adjoint and tangent linear system for general
problems, then consider the special case of mixed PDE-ODE systems and
finally discuss the specific adjoint and tangent linear versions of
multistage and Rush-Larson schemes.

We begin by considering a general system of discretized equations in the form:
find $y \in Y \subset \mathbb{R}^M$ such that
\begin{equation}
  \label{eq:AuF}
  F(y) = 0,
\end{equation}
The adjoint equation of \eqref{eq:AuF} associated with a real-valued
functional of interest $J: Y \rightarrow \R$ is: find the adjoint
solution $\bar{y} \in \mathbb R^M$ such that
\begin{equation}
  \label{eq:adjoint}
  \frac{\partial F}{\partial y}^{\ast}(y) \bar{y} = \frac{\partial J}{\partial y}(y)
\end{equation}
where the superscript $\ast$ denotes the adjoint operator. The
associated tangent linear equation with respect to some auxiliary
control parameter $m$ is given by: find the tangent linear solution
$\dot{y}$ such that
\begin{equation}
  \label{eq:tangent_linear}
  \frac{\partial F}{\partial y} (y)\dot{y}
    = \frac{\partial F}{\partial m}(y)
\end{equation}

We now turn our attention to the case where $F$ is a coupled PDE-ODE
system. For a guiding example, we consider again the bidomain
equations~\eqref{eq:bidomain:odes}--\eqref{eq:bidomain:pdes} and,
without loss of generality restrict ourselves to a single iteration
($N = 1$).  For brevity, we denote the ODE operator
\eqref{eq:bidomain:odes} as $O$ and the PDE operator
\eqref{eq:bidomain:pdes} as $P$. We can then rewrite the scheme in the
form \eqref{eq:AuF}:
\begin{subequations}
    \begin{align}
        v_0 - v^0 &= 0, &&\text{(Set initial condition)}\\
        s_0 - s^0 &= 0, &&\text{(Set initial condition)}\\
        O(v^*, s^*; v^0, s^0) &= 0, &&\text{(Solve bidomain ODEs \eqref{eq:bidomain:odes})}\\
        P(v^{\dag}, u_e^1; v^*) &= 0, &&\text{(Solve bidomain PDEs \eqref{eq:bidomain:pdes})}\\
        O(v^1, s^1; v^{\dag}, s^*) &= 0.&&\text{(Solve bidomain ODEs for final state)}
    \end{align}
    \label{eq:bidomain:structure}
\end{subequations}
The first pair of equations in~\eqref{eq:bidomain:structure} set the
initial conditions, the following equations represents a collection of
nonlinear ODEs, the third equation a system of PDEs, and the last
equation again a collection of ODEs.

From \eqref{eq:adjoint}, the adjoint system for
\eqref{eq:bidomain:structure} derives as:
\begin{equation}
    \left(
    \begin{array}{cc|cc|cc|cc}
        \multicolumn{2}{c|}{\multirow{2}{*}{$I$}}  &   \multicolumn{2}{|c|}{\multirow{2}{*}{$\frac{\partial{O^*}}{\partial{(v^0, s^0)}}$}} & \multicolumn{2}{c|}{\multirow{2}{*}{$0$}}  &   \multicolumn{2}{c}{\multirow{2}{*}{$0$}}   \\
     \phantom{xx}&  &   &   &  &   &   \\
        \hline
        \multicolumn{2}{c|}{\multirow{2}{*}{$0$}}   &       \multicolumn{2}{c|}{\multirow{2}{*}{$\frac{\partial{O^*}}{\partial{(v^*, s^*)}}$}} & \multicolumn{2}{c|}{\multirow{2}{*}{$\frac{\partial{P^*}}{\partial{(v^0, s^0)}}$}}& \multicolumn{2}{c}{\multirow{2}{*}{$\frac{\partial{O^*}}{\partial{s^*}}$}}\\
          & & & &   &   &   &   \\
        \hline
        \multicolumn{2}{c|}{\multirow{2}{*}{$0$}} &\multicolumn{2}{c|}{\multirow{2}{*}{$0$}}  &       \multicolumn{2}{c|}{\multirow{2}{*}{$\frac{\partial{P^*}}{\partial{(v^{\dag}, u_e^1)}}$}} & \multicolumn{2}{c}{\multirow{2}{*}{$\frac{\partial{O^*}}{\partial{v^\dag}}$}} \\
        &   &   &   & & &   &   \\
        \hline
        \multicolumn{2}{c|}{\multirow{2}{*}{$0$}}  &\multicolumn{2}{c|}{\multirow{2}{*}{$0$}}  & \multicolumn{2}{c|}{\multirow{2}{*}{$0$}} &  \multicolumn{2}{c}{\multirow{2}{*}{$\frac{\partial{O^*}}{\partial{(v^1, s^1)}}$}}  \\
        &   &   &   &   &   & &
  \end{array}
    \right)
  \begin{pmatrix}
    \bar v^0 \\ \bar s^0 \\ \bar v^{\ast} \\ \bar s^{\ast} \\ \bar v^{\dag} \\ \bar u_e^{1} \\ \bar v^{1} \\ \bar s^{1} \\
  \end{pmatrix}
    =
  \begin{pmatrix}
      0 \\ \multirow{2}{*}{\vdots} \\ \\  \\  \\ \\ \frac{\partial J}{\partial v^1} \\ \frac{\partial J}{\partial s^1}
  \end{pmatrix}
    \label{eq:bidomain:adjoint}
\end{equation}
For brevity, the left hand side matrix combines $2\times2$ blocks, and
the functional of interest is assumed to only depend on the final
state variables $v^1, s^1$. The adjoint system
\eqref{eq:bidomain:adjoint} is a linear coupled PDE-ODE system with
upper-triangular block structure, which can efficiently by solved by
backwards substitution.  The last block row represents a linear
adjoint ODE system, preceded by a adjoint PDE system, and preceded by
another adjoint ODE system, and finalised by a variable assignment in
the first block row which represents the adjoint version of setting
the initial conditions.

The derivation of adjoint and tangent linear equations for finite
element discretizations of PDEs, such as the third block row in
\eqref{eq:bidomain:adjoint}, is well-established, see
e.g.~\citep{FarrellEtAl2013}. However, the automated derivation of
adjoint and tangent linear systems for the multistage and Rush-Larsen
discretizations of the generic ODE system~\eqref{eq:ivp}, such as the
second and fourth block row in \eqref{eq:bidomain:adjoint}, is less
so, and these derivations are presented below.

\subsection{Adjoints and tangent linearizations of multistage schemes}

Consider a multistage discretization of~\eqref{eq:ivp} as described in
Section~\ref{sec:multistage:schemes} with a given Butcher tableau
$a_{ij}, b_j, c_i$, $i, j = 1, \dots, s$. Taking $s = 3$ for
illustrative purposes, we can
write~\eqref{eq:multistage:stage}--\eqref{eq:multistage:final} in form
\eqref{eq:AuF} as:
\begin{equation*}
  \begin{pmatrix}
    y^0 \\
    k_1 \\
    k_2 \\
    k_3 \\
    y^{1}
  \end{pmatrix}
  -
  \begin{pmatrix}
    y_0 \\
    f(y^0 + \kappa w_1, T_0 + c_1 \kappa) \\
    f(y^0 + \kappa w_2, T_0 + c_2 \kappa) \\
    f(y^0 + \kappa w_3, T_0 + c_3 \kappa) \\
      y^0 + \kappa \sum_{i=1}^3 b_i k_i
  \end{pmatrix}
    =0,
\end{equation*}
with
\begin{equation*}
  w_i = y^0 + \kappa \sum_{j=1}^s a_{ij} k_j.
\end{equation*}
Assuming that the functional of interest $J$ is independent of
the internal stage values $k_i$, we can derive the adjoint problem as:
\begin{equation*}
  \begin{pmatrix}
    \bar{y}^0 \\
    \bar{k}_1 \\
    \bar{k}_2 \\
    \bar{k}_3 \\
    \bar{y}^{1}
  \end{pmatrix}
  -
  \begin{pmatrix}
    0 & \frac{\partial f_1}{\partial w_1}^* & \frac{\partial f_2}{\partial w_2}^* & \frac{\partial f_3}{\partial w_3}^* & I \\
    0 &  \kappa a_{11} \frac{\partial f_1}{\partial w_1}^* & \kappa a_{21} \frac{\partial f_2}{\partial w_2}^* & \kappa a_{31} \frac{\partial f_3}{\partial w_3}^* & \kappa b_1 \\
    0 & 0                                                 &   \kappa a_{22} \frac{\partial f_2}{\partial w_2}^*& \kappa a_{32} \frac{\partial f_3}{\partial w_3}^* & \kappa b_2 \\
    0 & 0                                                 &  0                                                 &  \kappa a_{33} \frac{\partial f_3}{\partial w_3}^* & \kappa b_3 \\
    0 & 0 & 0 & 0 & 0 \\
  \end{pmatrix}
  \begin{pmatrix}
    \bar{y}^0 \\
    \bar{k}_1 \\
    \bar{k}_2 \\
    \bar{k}_3 \\
    \bar{y}^{1}
  \end{pmatrix}
  =
  \begin{pmatrix}
    \frac{\partial J}{\partial y^{0}} \\
    0 \\
    0 \\
    0 \\
    \bar{y}_0
  \end{pmatrix}
  ,
\end{equation*}
where $\bar{y}_0$ is the terminal condition for the adjoint solution
at $T_1$.

Generalizing to $s$ stages, we see that we first solve for the adjoint
stage values and then compute the adjoint solution at time $T_0$ via
\begin{subequations}
  \label{eq:multistage:adjoint}
  \begin{align}
\left(I - \kappa a_{ii} \frac{\partial f_i}{\partial w_i}^*\right) \bar{k}_i &= \kappa b_i \bar{y}^{1} + \sum_{\substack{j=i + 1}}^s \kappa a_{ji} \frac{\partial f_j}{\partial w_j}^* \bar{k}_j, \\
\bar{y}^0 &= \bar{y}^{1} + \sum_{i=1}^s \frac{\partial f_i}{\partial w_i}^* \bar{k}_i.
\end{align}
\end{subequations}

Following a similar calculation, the tangent linearization of the
multistage scheme with $s$ stages is: given $y_0$ at $T_0$ compute
$\dot{k}_i$ for $i = 1, \dots, s$ and then $\dot{y}^{1}$ at $T_1$ via
\begin{subequations}
  \label{eq:multistage:tangent_linear}
  \begin{align}
      \left(I - \kappa a_{ii} \frac{\partial f_i}{\partial w_i}\right) \dot{k}_i & = \frac{\partial f_i}{\partial w_i} \dot{y}^n + \sum_{\substack{j = 1}}^{i-1} \kappa a_{ij} \frac{\partial f_i}{\partial w_i} \dot{k}_j + \frac{\partial f_i}{\partial m}, \\
    \dot{y}^{n+1} &= \dot{y}^n + \sum_{i=1}^s \kappa b_i \dot{k}_i.
  \end{align}
\end{subequations}

\subsection{Adjoints and tangent linearizations of Rush-Larsen schemes}

We now turn to consider adjoints and tangent linearizations of the
(generalized) Rush-Larsen schemes. We here derive the equations for
the first order generalized Rush-Larsen scheme given
by~\eqref{eq:grl1}:
\begin{equation*}
  y_i^{1} = y_i^{0} + J_{i}^{-1} f_i^0 \left(e^{\kappa J_{i} }-1 \right)
\end{equation*}
for each component $i = 1, \dots, M$ of the state variable.

We can write~\eqref{eq:grl1} in form \eqref{eq:AuF} as
\begin{equation*}
\begin{pmatrix}
  y^0 - y_0 \\
  y^1 - L(y^0) - y^0 \\
\end{pmatrix}
= 0
\end{equation*}
where $y_0$ is the given initial condition at $T_0$ and $L(y^0) = \{ L_i(y^0) \}_{i=1}^M$ with
\begin{equation*}
  L_i(y^0)
  = J_{i}^{-1} f_i^0 \left( e^{\kappa J_{i} }-1 \right )
  \equiv \frac{\partial f_i}{\partial y_i}(y^0, T_0)^{-1} f_i(y_0, T_0) (e^{\kappa \frac{\partial f_i}{\partial y_i}(y^0, T_0)} - 1) .
\end{equation*}
The adjoint system of the first order generalised Rush-Larsen scheme
is then given by
\begin{equation}
  \label{eq:grl1:adjoint}
  \begin{pmatrix}
  \bar{y}^0 \\
  \bar{y}^1 \\
\end{pmatrix}
    -
  \begin{pmatrix}
    0 & I +  \frac{\partial L}{\partial y^0}^{\ast}\\
    0 & 0 \\
  \end{pmatrix}
  \begin{pmatrix}
  \bar{y}^0 \\
  \bar{y}^1 \\
\end{pmatrix}
=
\begin{pmatrix}
  \frac{\partial J}{\partial y^0} \\
  \bar{y}_0 \\
\end{pmatrix}
\end{equation}
where $\bar{y}_0$ is a given terminal condition for the adjoint at
$T_1$. Similarly, the tangent linear system \eqref{eq:tangent_linear}
with respect to an auxiliary parameter $m$ is given by
\begin{equation}
  \label{eq:grl1:tangent_linear}
  \begin{pmatrix}
    \dot{y}^0 \\
    \dot{y}^1
  \end{pmatrix}
    -
  \begin{pmatrix}
    0 & 0 \\
    I + \frac{\partial L}{ \partial y^0} & 0 \\
  \end{pmatrix}
  \begin{pmatrix}
    \dot{y}^0 \\
    \dot{y}^1
  \end{pmatrix}
  =
  \begin{pmatrix}
    \dot{y}_0 \\
    \frac{\partial L(y^0)}{\partial m}
  \end{pmatrix}
\end{equation}
where $\dot{y}_0$ is a given initial condition at $T_0$.

The adjoint and tangent linear equations for the other Rush-Larsen
schemes are derived following the same steps, and we therefore omit
the details here.

\section{FEniCS and Dolfin-adjoint abstractions and algorithms}
\label{sec:implementation}

The FEniCS Project defines a collection of software components
targeting the automated solution of differential equations via finite
element methods~\citep{LoggMardalEtAl2012a}. The components 
include the Unified Form Language (UFL)~\citep{AlnaesEtAl2014}, the
FEniCS Form Compiler (FFC)~\citep{LoggOlgaardEtAl2012} and the finite
element library DOLFIN~\citep{LoggWells2010}. The separate
dolfin-adjoint project and software automatically derives the discrete
adjoint and tangent linear models from a forward model written in the
Python interface to DOLFIN~\citep{FarrellEtAl2013}. This section
presents our extensions of the FEniCS form language and the FEniCS
form compiler, and other new FEniCS and dolfin-adjoint software
features targeting coupled PDE-ODE systems in general and collections
of systems of ODE in particular.

\subsection{Variational formulation of collections of ODE systems}

Consider a spatial domain $\Omega \subset \R^d$ tessellated by a mesh
$\triang_h$, and a collection of general initial problems as defined
by~\eqref{eq:intro:ivp} over a set of points $X$ defined relative to
$\triang_h$. For $x_i \in X$, denote by $\delta_{x_i}$ the Dirac delta
function centered at $x_i \in \R^d$ such that
\begin{equation*}
  f(x_i)
  = \int_{\R^d} f \, \delta_{x_i} \dx
\end{equation*}
for all continuous, compactly supported functions $f$. We will also
write
\begin{equation}
  \sum_{x_i \in X} \int_{\R^d} f \, \delta_{x_i} \dx
  \equiv \sum_{x_i \in X} \int f \dP(x_i)
  \equiv \sum_{x_i \in X} \inner{f}{1}_{x_i}
  \equiv \inner{f}{1}_X
  \label{eq:point_cloud_integration}
\end{equation}

Drawing inspiration from our context of finite element variational
formulations, we can then write~\eqref{eq:intro:ivp} as: find
$y(\cdot, t)$ such that
\begin{equation}
  \label{eq:ivp:variational}
  \inner{y_t(t)}{\psi_i}_{X} = \inner{f(y, t)}{\psi_i}_{X}
\end{equation}
for all (basis) functions (or distributions) $\psi_i$ such that
$\psi_{i}(x_i) = 1$ and $\psi_{j}(x_i) = 0$ for $j \not = i$.

The remainder of this section describes the extensions of the FEniCS
and Dolfin-adjoint systems to allow for in particular abstract
representation and efficient forward and reverse solution of systems
of the form~\eqref{eq:ivp:variational}.

\subsection{Extending UFL with vertex integrals}

UFL is an expressive domain-specific language for abstractly
representing (finite element) variational formulations of differential
equations. In particular, the language defines syntax for integration
over various domains. Consider a mesh $\triang_h$ of geometric
dimension $d$ with cells $\{\mathcal{T}\}$, interior facets
$\{\mathcal{F}^i\}$ and boundary facets $\{\mathcal{F}\}^b$. Interior
facets are defined as the $d-1$ dimensional intersections between two
cells and thus for any interior facet $\mathcal{F}^i$ we can write
$\mathcal{F}^i = \mathcal{T}^{+} \cap \mathcal{T}^{-}$. UFL defines
the sum of integrals over cells, sum of integrals over boundary facets
and sum of integrals over interior facets by the \texttt{dx}, \texttt{ds},
and \texttt{dS} measures, respectively. For example, the following linear
variational form defined in terms of a piecewise polynomial $u$
defined relative over $\triang_h$:
\begin{equation*}
  L(u) = \sum_{\{\mathcal{T}\}} \int_{\mathcal{T}} u  \dx
  + \sum_{\{\mathcal{F}^i\}} \int_{\mathcal{F}^i} u |_{T^{+}}  \ds
  + \sum_{\{\mathcal{F}^b\}} \int_{\mathcal{F}^b} u \ds
\end{equation*}
is naturally expressed in UFL as:
\begin{python}
L = u*dx + u('+')*dS + u*ds
\end{python}

To allow for point evaluation over the vertices of a mesh
(cf.~\eqref{eq:point_cloud_integration}), we have introduced a new
\emph{vertex integral} (or \emph{vertex measure} in UFL terms) type
with default instantiation \texttt{dP}:
\begin{python}
from ufl import Measure
dP = Measure('vertex')
L = f*dP
\end{python}
In agreement with~\eqref{eq:point_cloud_integration}, the vertex
integral is defined by:
\begin{equation}
  \mathtt{f*dP} \equiv \inner{f}{1}_{\mathcal{V}(\triang_h)} = \sum_{x_i \in \mathcal{V}(\triang_h)} f (x_i),
\end{equation}
where $V(\triang_h)$ is the set of all vertices of the mesh
$\triang_h$. Vertex measures restricted to a subset of vertices can be
defined as for all other UFL measure types.

This basic extension of UFL crucially allows for the abstract
specification of collections of ODEs such
as~\eqref{eq:ivp:variational} in a manner that is consistent and
compatible with abstract specification of finite element variational
formulations of PDEs. It also allows for the native specification of
e.g. point sources in finite element formulations of PDEs in FEniCS.

\subsection{Vertex integrals in the FEniCS Form Compiler}

The FEniCS Form Compiler FFC generates specialized C++ code~\citep{AlnaesLoggEtAl2009} from the symbolic UFL
representation of variational forms and finite element
spaces~\citep{LoggOlgaardEtAl2012}. To accommodate the new vertex
integral type, we have extended the UFC interface with a class
\texttt{vertex\_integral} defining the interface for the tabulation of
the element matrix corresponding to the evaluation of an expression at
a given vertex, listed below:
\begin{cpp}
/// This class defines the interface for the tabulation of
/// an expression evaluated at exactly one point.
class vertex_integral: public integral
{
public:

  /// Constructor
  vertex_integral() {}

  /// Destructor
  virtual ~vertex_integral() {}

  /// Tabulate tensor for contribution from local vertex
  virtual void tabulate_tensor(double * A,
                               const double * const * w,
                               const double * coordinate_dofs,
                               std::size_t vertex,
                               int cell_orientation) const = 0;
};
\end{cpp}
The FFC code generation pipeline has been correspondingly extended to
allow for the generation of optimized code from the UFL
representation of variational forms involving vertex integrals. This
allows for subsequent automated assembly of variational forms
involving single vertex integrals (as illustrated above) and also in
combination with other (cell, interior facet, exterior facet)
integrals.

\subsection{DOLFIN features for solving collections of ODE systems}

\subsubsection{Assembly of vertex integrals}

We have extended DOLFIN with support for automated assembly of variational forms that include vertex
integrals. The support is currently limited to forms
defined over test and trial spaces with vertex-based degrees of
freedom only. The assembly algorithm follows the standard finite
element assembly pattern by iterating over the vertices of the mesh,
computing the map from local to global degrees of freedom, evaluating
the local element tensor based on generated code,
and adding the contributions to the global tensor. The vertex integral
assembler runs natively in parallel via MPI.

\subsubsection{Specification and generation of multistage schemes}

Since version \version, DOLFIN supports the specification of multistage
schemes of the
form~\eqref{eq:multistage:stage}--\eqref{eq:multistage:final} for the
solution of collections of ODE systems of the
form~\eqref{eq:intro:ivp}, via their Butcher tableaus. The DOLFIN
class \texttt{ButcherMultiStageScheme} takes as input the right hand
side expression $f$, the solution $y$, the Butcher tableau specified
via $a$, $b$ and $c$, and secondary variables such as a time variable,
the (integer) order of the scheme. From this specification, DOLFIN automatically generates a
variational formulation of each of the separate stages in the
multistage scheme, with each stage in accordance
with~\eqref{eq:multistage:stage}-\eqref{eq:multistage:final}. DOLFIN
can also automatically generate the variational formulations
corresponding to the adjoint scheme~\eqref{eq:multistage:adjoint} and
to the tangent linear scheme~\eqref{eq:multistage:tangent_linear} on
demand. A set of common multistage schemes are predefined including
Crank-Nicolson, explicit Euler, implicit Euler, 4$^{\mathrm{th}}$-order
explicit Runge-Kutta, ESDIRK3, ESDIRK4~\citep{Butcher2008,Kvaerno2004}.

Similar features have also been implemented for easy specification of
Rush-Larsen schemes cf.~\eqref{eq:rl1}--\eqref{eq:grl2} and the
automated generation of the corresponding adjoint and tangent linear
schemes cf.~e.g.~\eqref{eq:grl1:adjoint} and
\eqref{eq:grl1:tangent_linear} through the class
\texttt{RushLarsenScheme}. Both
\texttt{ButcherMultiStageScheme} and \texttt{RushLarsenScheme}
subclass the \texttt{MultiStageScheme} class.

\subsubsection{PointIntegralSolvers}

To allow for the efficient solution of the multistage schemes for
collections of systems of ODEs, we have introduced a targeted
\texttt{PointIntegralSolver} class in DOLFIN \version. This solver
class takes a \texttt{MultiStageScheme} as input and its main
functionality is to compute the solutions over a single time step.

The solver iterates over all vertices of the mesh and solves for the
relevant (stage and/or final) variables at each vertex. The solution
algorithm for each stage depends on whether the stage is explicit or
implicit. For implicit stages, a custom Newton solver is invoked that
allows for Jacobian reuse across stages, across vertices, and/or
across time steps on demand. As the resulting linear systems are
typically small and dense, a direct (LU) algorithm is used for the
inner solves. For explicit stages, a simple vector update is
performed.

The point integral solver runs natively in parallel via MPI. For a
cell-partitioned mesh distributed between $N$ processes, each cell is
owned by one process and that process performs the solve for each
vertex on the cell (once). As little communication is required
between processes, the total solve is expected to scale
linearly in $N$.

\subsection{Extensions to dolfin-adjoint}

Dolfin-adjoint has been extended to support the new features including
automatically deriving and computing adjoint and tangent linear
solutions for \texttt{PointIntervalSolver}. The symbolic derivation of
the variational formulation for the adjoint and tangent linear
equations for the \texttt{MultiStageScheme} are based
on~\eqref{eq:multistage:adjoint},~\eqref{eq:multistage:tangent_linear}
for the Butcher tableau defined schemes
and~\eqref{eq:grl1:adjoint},~\eqref{eq:grl1:tangent_linear} and its
analogies for the Rush-Larsen schemes.

The overall result is that operator splitting algorithms as described
in Section~\ref{sec:operator_splitting} can be fully specified and
efficiently solved within FEniCS. Moreover, the corresponding adjoint
and tangent linear models and first and second order functional
derivatives may be efficiently computed using dolfin-adjoint with only
a few lines of additional code.

\section{Applications}
\label{sec:applications}

In this section, we demonstrate the applicability and performance of
our implementation by considering the two examples presented in the
introduction, originating from the computational modelling of
mitochondrial swelling and cardiac electrophysiology respectively. The
complete supplementary code is openly available, see~\citep{TheCode}.

\subsection{Application: mitochondrial swelling}

We consider the mathematical model defined by~\eqref{eq:mitochondria}.
Inspired by~\citep{Eisenhofer2013}, we
let $x = (x_0, x_1) \in \Omega = [0, 1]^2$ and $t \in [0, T]$ with $T
= 35$, and consider the initial conditions $N_{1, 0} = 1$, $N_{2,
  0} = 0$, $N_{3, 0} = 0$ and
\begin{equation}
  u_0(x) = \frac{30}{\int_{\Omega} M(x) \dx} M(x), \quad M(x_0, x_1) =
  \frac{1}{2 \pi} e^{- 0.5 (X_0^2 + X_1^2)},
\end{equation}
with $X_0 = (\alpha - \beta) x_0 + \beta$, $X_1 = (\alpha - \beta) x_1
+ \beta$ taking $\alpha = 3$ and $\beta = -1$. We consider
\begin{equation}
  f(s) =
  \begin{cases}
    0 & s < C^{-}, \\
    f^{\ast} & s > C^{+}, \\
    \frac{f^{\ast}}{2} \left (1 -  \cos \left (\frac{s - C^{-}}{C^{+} - C^{-}} \pi \right ) \right ) & \text{otherwise},
  \end{cases}
\end{equation}
and
\begin{equation}
  g(s) =
  \begin{cases}
    g^{\ast} & s > C^{+}, \\
    \frac{g^{\ast}}{2} \left ( 1 - \cos \left ( \frac{s}{C^{+}} \pi \right ) \right )  & \text{otherwise}.
  \end{cases}
\end{equation}
with $C^{-} = 20$, $C^{+} = 200$, $f^{\ast} = 1$, $g^{\ast} = 0.1$,
$d_1 = 2 \times 10^{-6}$, and $d_2 = 30$.

We are interested in the total amount of completely swollen
mitochondria (whose density is given by $N_3$) at the final time $T$ and its
sensitivity to the initial condition for $u$ ($u_0$). Thus, our
functional of interest $J$ is given by
\begin{equation}
  \label{eq:mitochondria:J}
  J = \int_{\Omega} N_3(T)
\end{equation}
and our goal is to compute the derivative of $J$ with respect to $u_0$.

We use a second order Strang splitting scheme
for~\eqref{eq:mitochondria}, a Crank-Nicolson discretization
in time and continuous piecewise linear finite elements in space.
The scheme then
reads as: given initial conditions $u_0$, $N_{1,0}$, $N_{2,0}$,
$N_{3,0}$, and time points $\{t_0, t_1, \dots, t_N\}$ with time step
$\kappa_n = t_{n+1} - t_n$, then for each $n = 0, 1, \dots, N - 1$:
\begin{compactenum}
\item
  Compute $u^\ast$ and $N_{i, 0}^\ast$ $i = 1, 2, 3$ solving
  \begin{equation}
    \label{eq:mitochondria:odes}
      \ddt{u} = d_2 g(u) N_2, \;
      \ddt{N_1} = - f(u) N_1, \;
      \ddt{N_2} = f(u) N_1 - g(u) N_2, \;
      \ddt{N_3} = g(u) N_2,
  \end{equation}
  over $\Omega \times [t_n, t_n + \frac{1}{2} \kappa_n]$ with initial
  conditions $u^{n}, N_{i}^{n}$ for $i = 1, 2, 3$.
\item
  Compute a solution $u^{\dagger}$ of: find $u_h \in V_h$ such that
  \begin{equation}
    \label{eq:mitochondria:pde}
    \int_{\Omega} \left ( u_h - u^{\ast} \right ) v
     + \kappa_n  \Grad \mathcal{A} (\avg{u_h})  \cdot \Grad v \dx = 0
  \end{equation}
  for all $v \in V_h$ with $\avg{u} = \frac{1}{2} \left (u + u^{\ast}
  \right )$ and where $\mathcal{A}(u) = u$ or $\mathcal{A}(u) =
  |u|^{q-2} u$.
\item
  Compute solutions $u^{n+1}$ and $N_{i}^{n+1}$, $i = 1, 2, 3$,
  solving~\eqref{eq:mitochondria:odes} over $\Omega \times [t_n +
    \frac{1}{2} \kappa_n, t_{n+1}]$ with initial conditions
  $u^{\dagger}$ and $N_{i}^{\ast}$, $i = 1, 2, 3$.
\end{compactenum}
We choose to discretize~\eqref{eq:mitochondria:odes} via the ESDIRK4
method in time and set a tolerance of $10^{-10}$ for the inner Newton
solves. Further, we take $q = 3$, take $\kappa_n = 0.5$, and let
$\triang_h$ be a uniform tessellation of $\Omega$ with $N_x \times N_x
\times 2$ triangles.

\begin{figure}
  \centering
  \begin{subfigure}[t]{0.24\textwidth}
    \centering
    \includegraphics[width=\textwidth]{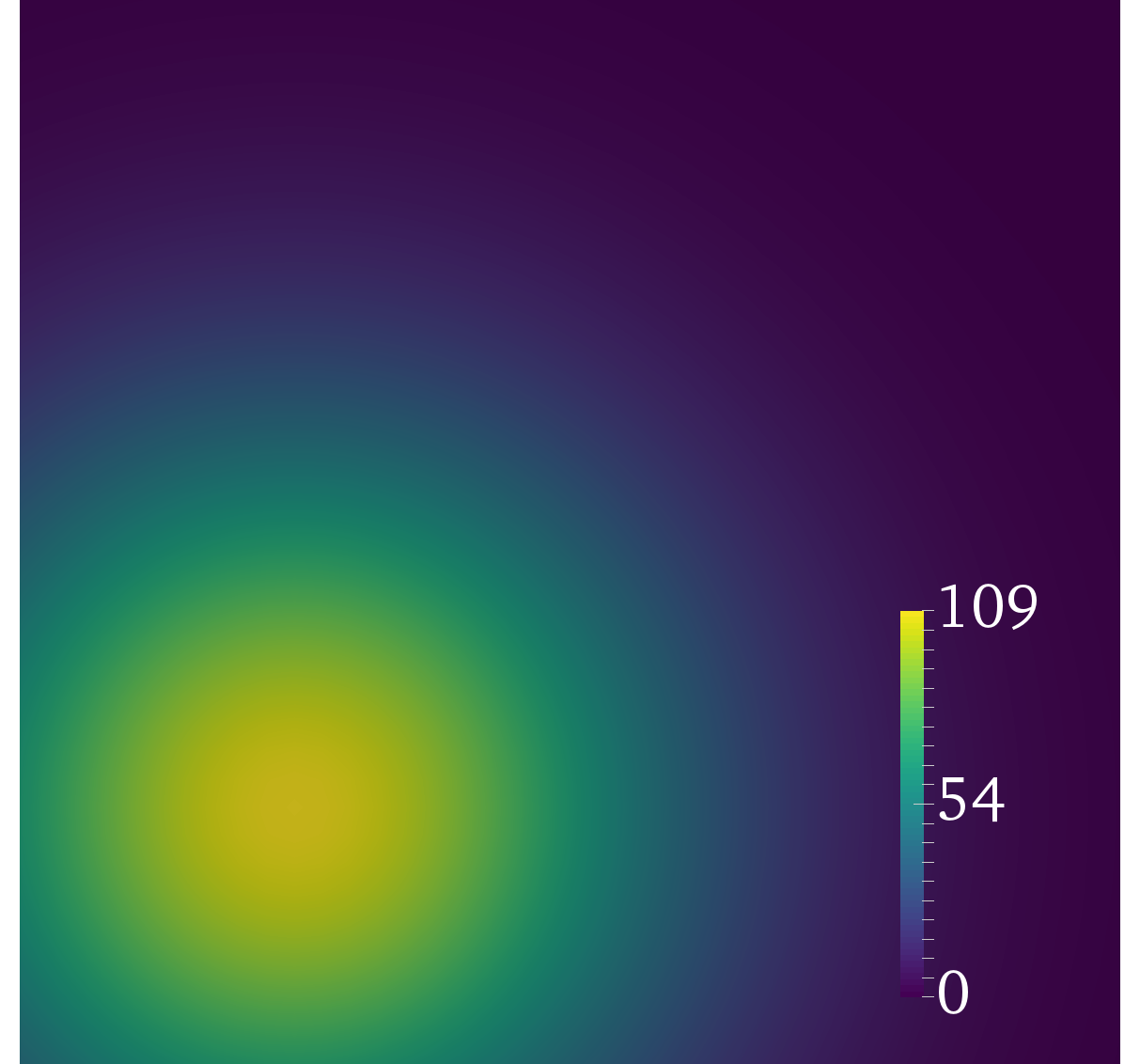}
    \caption{$u_0$}
  \end{subfigure}
  \begin{subfigure}[t]{0.24\textwidth}
    \centering
    \includegraphics[width=\textwidth]{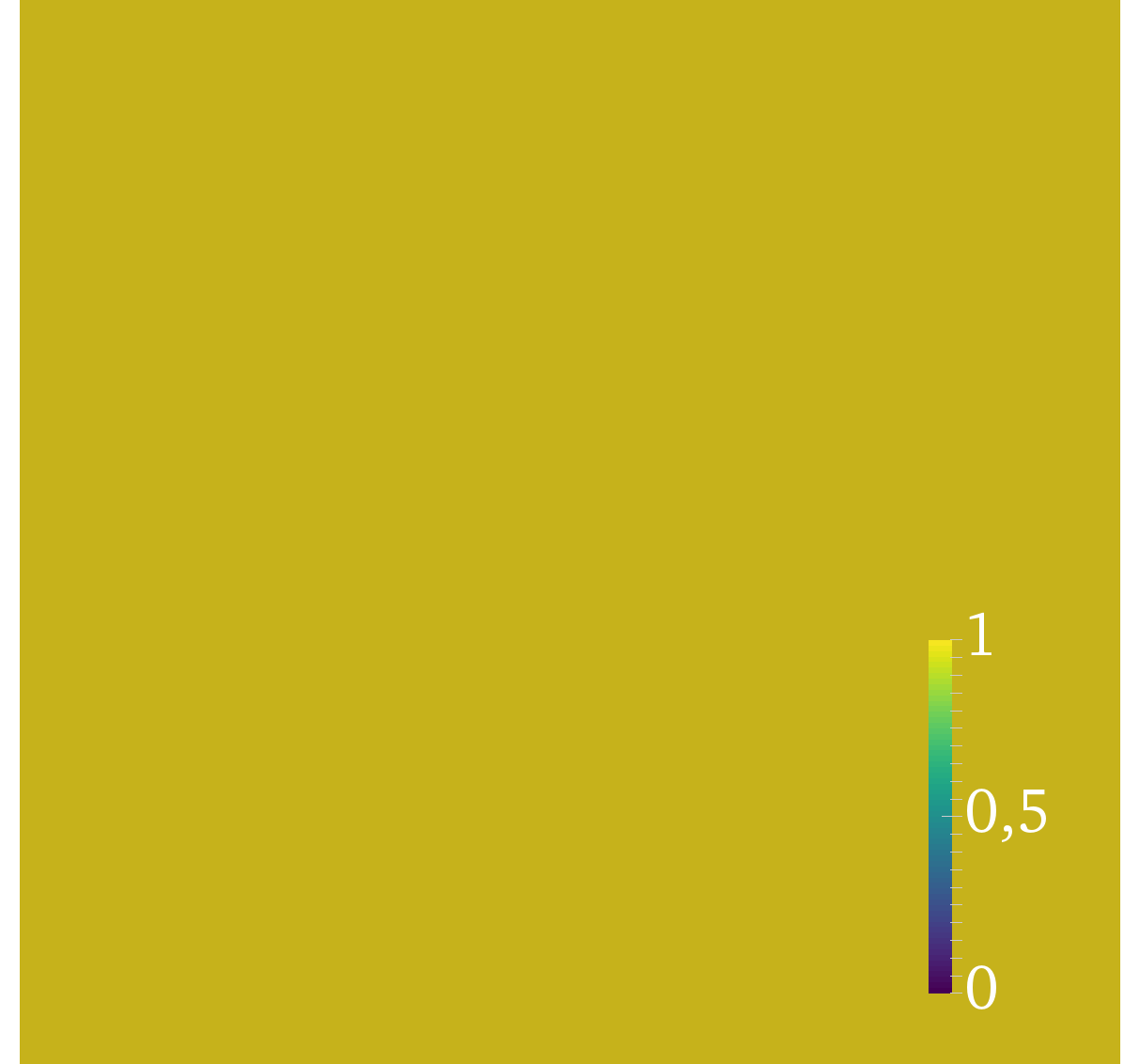}
    \caption{$N_{1, 0}$}
  \end{subfigure}
  \begin{subfigure}[t]{0.24\textwidth}
    \centering
    \includegraphics[width=\textwidth]{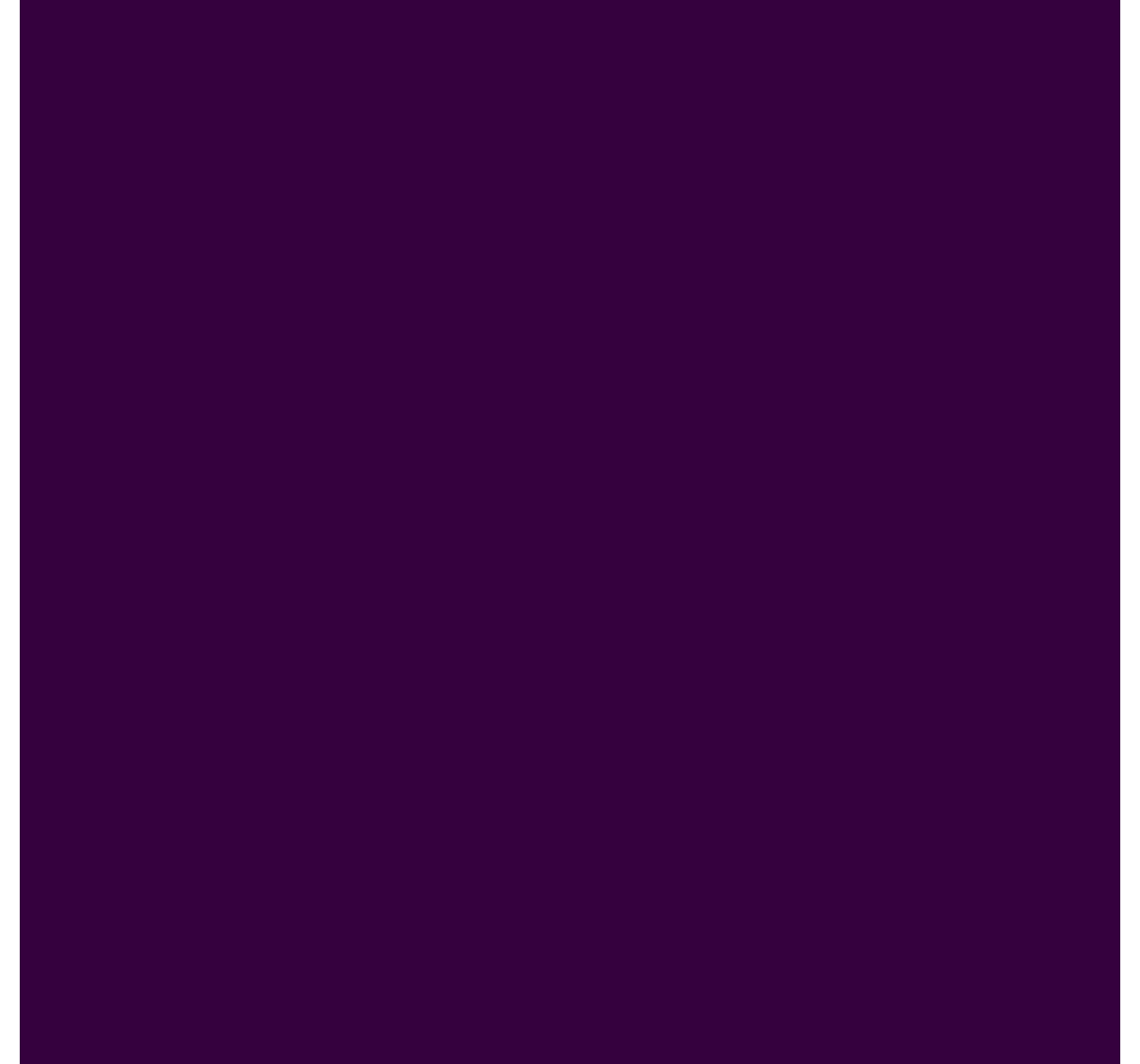}
    \caption{$N_{2, 0}$}
  \end{subfigure}
  \begin{subfigure}[t]{0.24\textwidth}
    \centering
    \includegraphics[width=\textwidth]{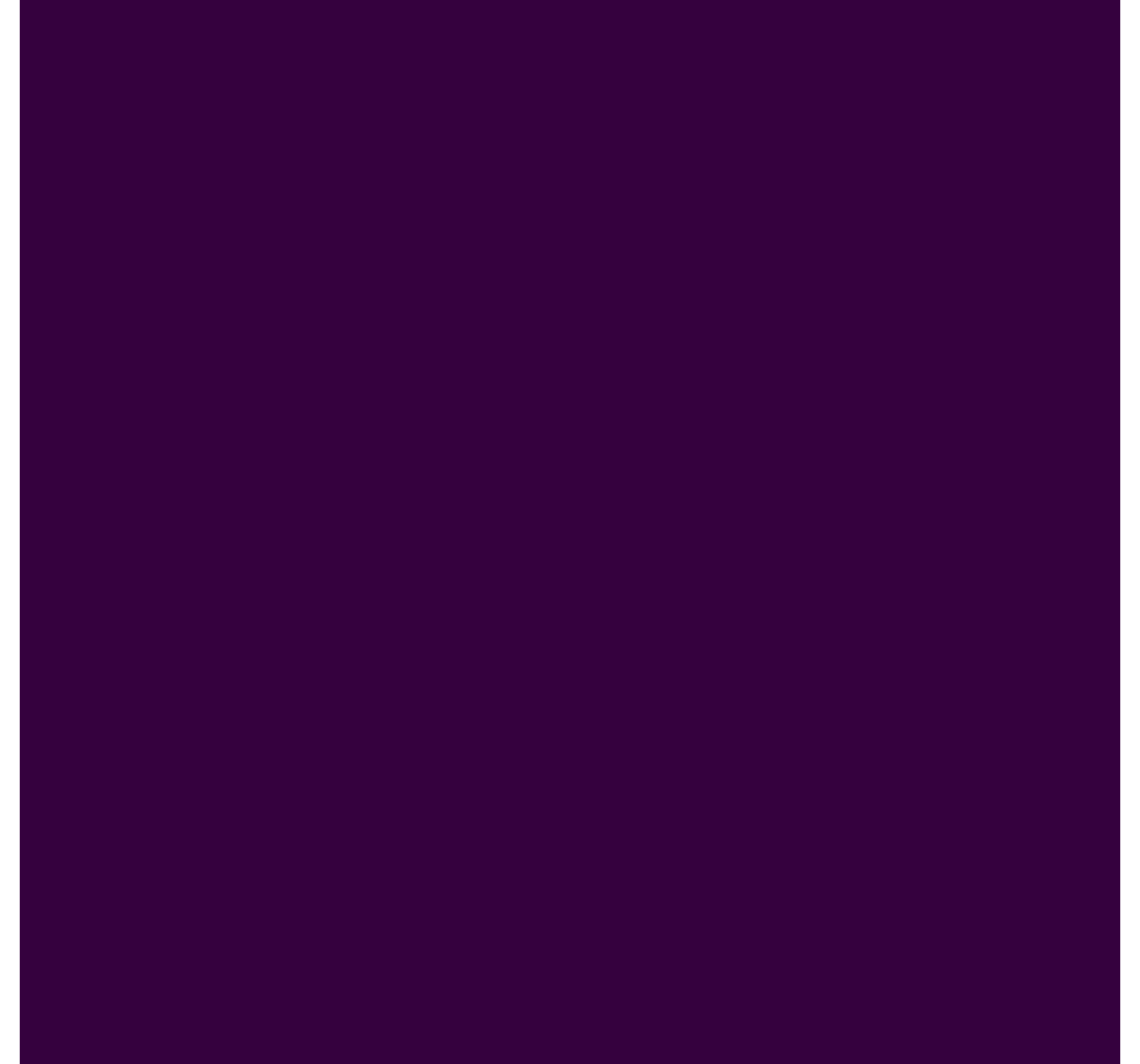}
    \caption{$N_{3, 0}$}
  \end{subfigure}
  \begin{subfigure}[t]{0.24\textwidth}
    \centering
    \includegraphics[width=\textwidth]{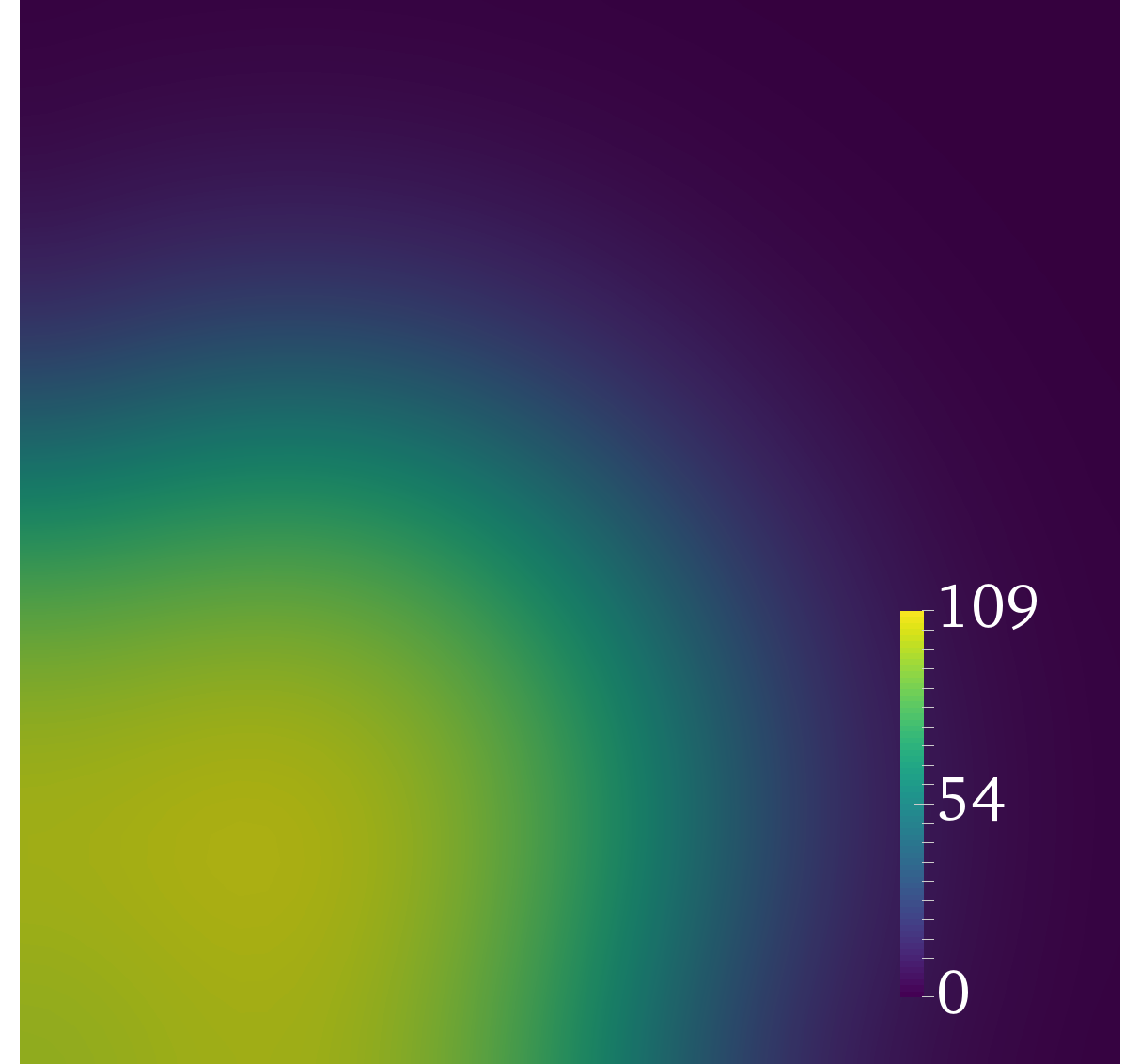}
    \caption{$u(T)$}
  \end{subfigure}
  \begin{subfigure}[t]{0.24\textwidth}
    \centering
    \includegraphics[width=\textwidth]{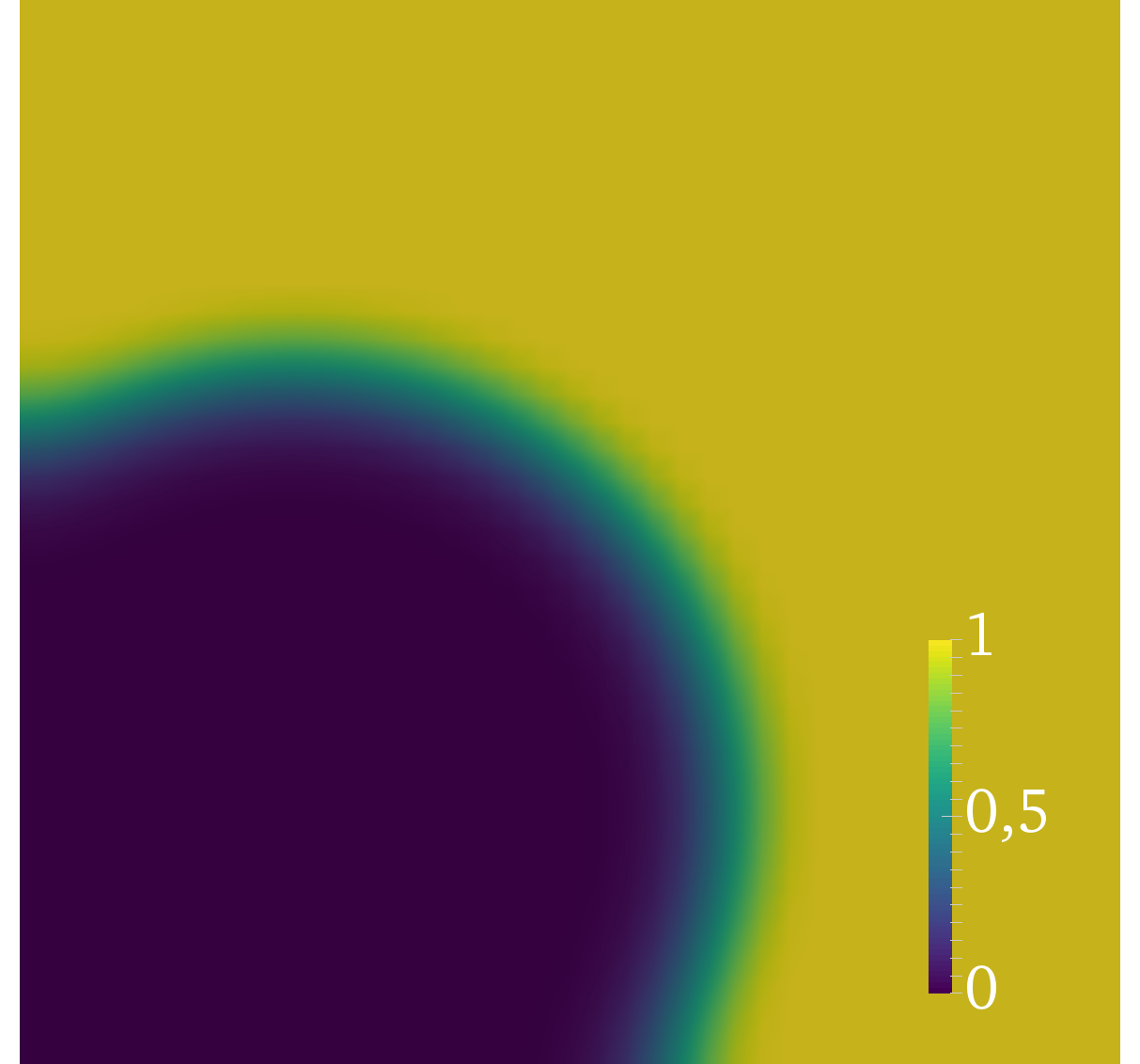}
    \caption{$N_{1}(T)$}
  \end{subfigure}
  \begin{subfigure}[t]{0.24\textwidth}
    \centering
    \includegraphics[width=\textwidth]{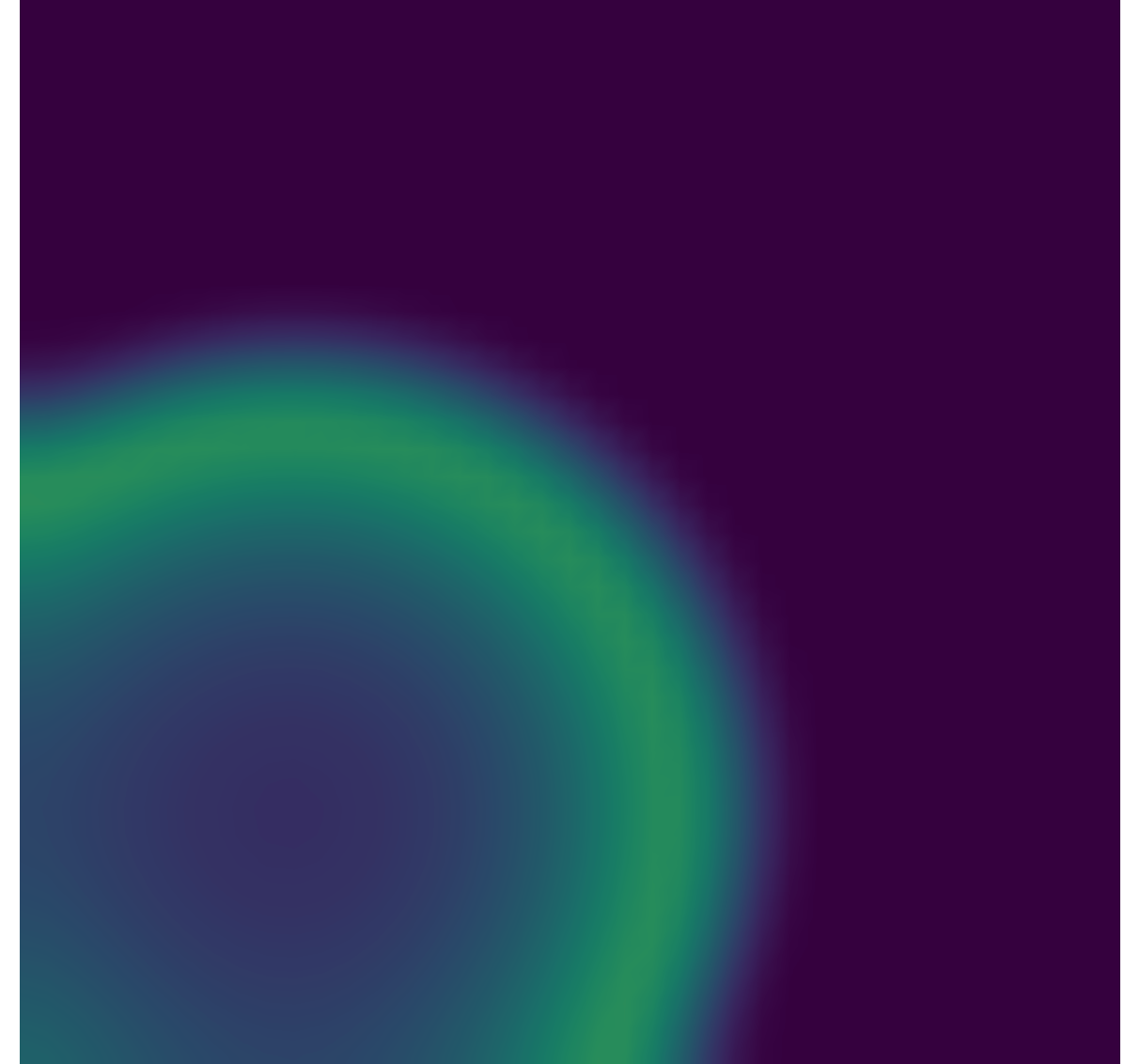}
    \caption{$N_{2}(T)$}
  \end{subfigure}
  \begin{subfigure}[t]{0.24\textwidth}
    \centering
    \includegraphics[width=\textwidth]{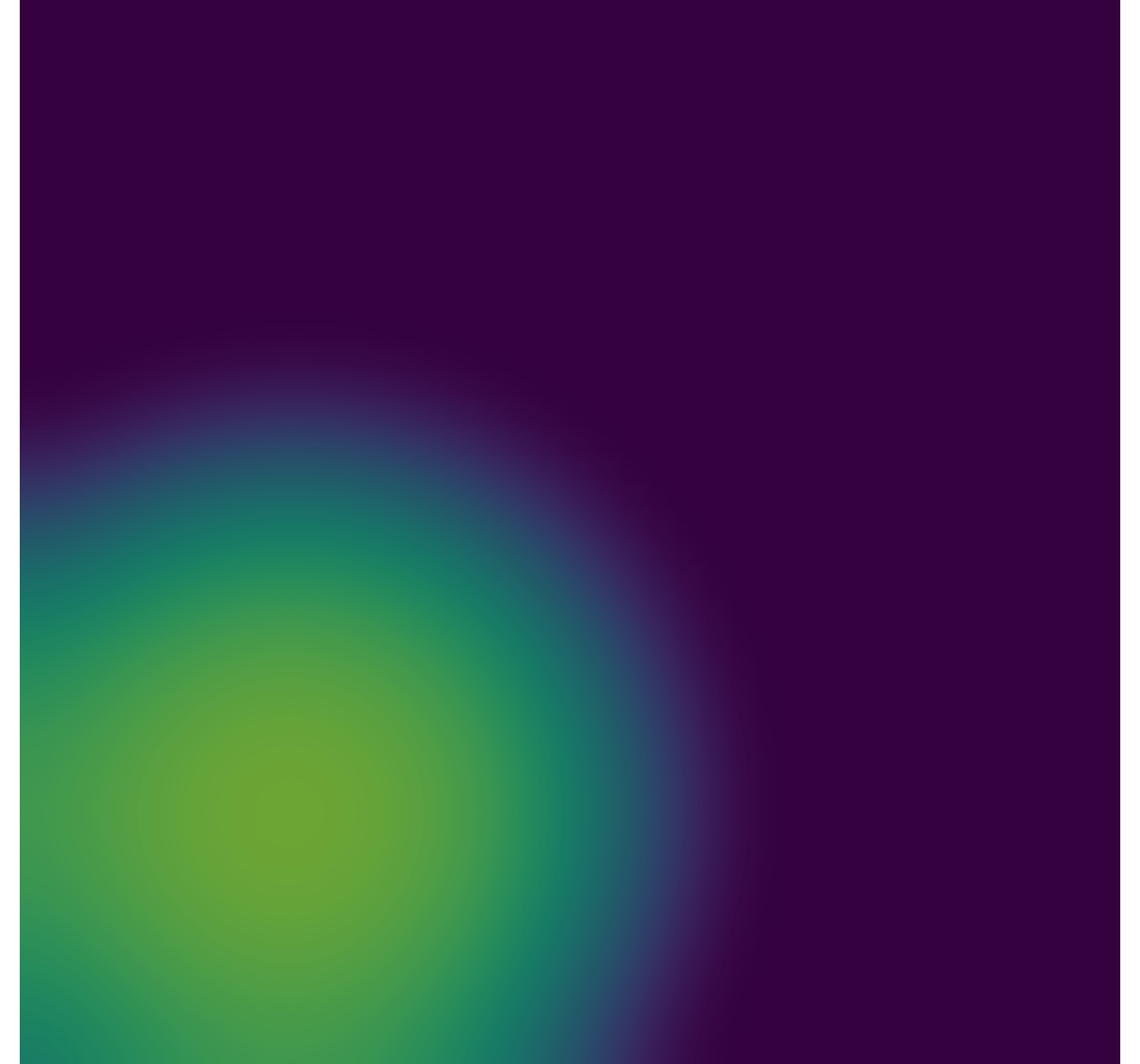}
    \caption{$N_{3(T)}$}
  \end{subfigure}
  \caption{Initial conditions and solutions $t = T = 35$ for the
    mitochondrial swelling example with $N_x = 40$. The scale is
    common for A and E and ranges from $0$ (blue) to 108.3
    (yellow). The scale is common for B-D and F-G and ranges from $0$
    (blue) to $1$ (yellow).}
  \label{fig:mito:solutions}
\end{figure}

\begin{figure}
  \centering
  \includegraphics[width=0.45\textwidth]{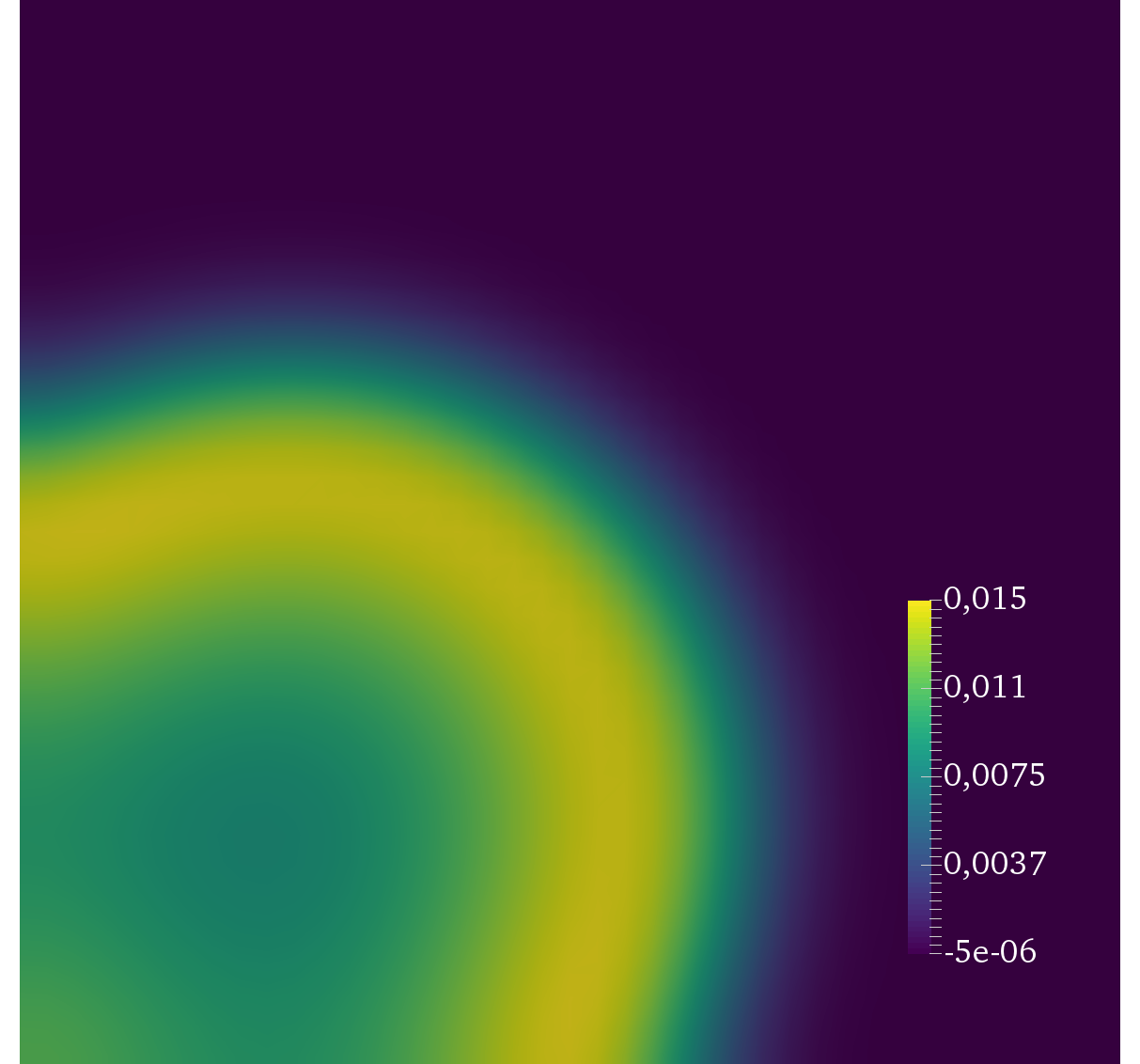}
  \includegraphics[width=0.45\textwidth]{u0.png}
  \caption{Sensitivity of $J$, the total amount of completely swollen
    mitochondria at $T$, $L^2$-gradient of $J$ (left) with respect to
    the initial condition $u_0$ (right).}
  \label{fig:mito:sensitivity}
\end{figure}

The initial and final solutions are presented in
Figure~\ref{fig:mito:solutions} while the $L^2$-gradient of the
objective functional $J$ given by~\eqref{eq:mitochondria:J} with
respect to the initial calcium concentration $u_0$ is presented in
Figure~\ref{fig:mito:sensitivity}. 

To verify that the computed gradient is correct, we performed a Taylor
test at a given spatial and temporal resolution with a given
perturbation seed. The results are listed in
Table~\ref{tab:mito:taylor} and demonstrate the expected orders of
convergence, indicating a correctly computed gradient.
\begin{table}
  \centering
  \begin{tabular}{l|cccc}
    \toprule
    $\delta u_0$ & $R_0 (\delta u_0)$ & order & $R_1 (\delta u_0)$ & order \\
    \midrule
    0.5 & $1.22 \times 10^{-3}$  &      & $1.25 \times 10^{-6}$ & \\
    0.25 & $6.14 \times 10^{-4}$   & 1.00 & $3.11 \times 10^{-7}$ & 2.00 \\
    0.125 & $3.07 \times 10^{-4}$   & 1.00 & $7.77 \times 10^{-8}$ & 2.00 \\
    0.0625 & $1.53 \times 10^{-4}$  & 1.00 & $1.94 \times 10^{-8}$ & 2.00 \\
    0.03125 & $7.67 \times 10^{-5}$ & 1.00 & $4.85 \times 10^{-9}$ & 2.00 \\
    \bottomrule
  \end{tabular}
  \vspace{1em}
  \caption{Taylor reminders $R_0 = |J(u_0 + \delta u_0)|$ and $R_1 =
    |J(u_0 + \delta u_0) - J(u_0) - \nabla J(u_0) \delta u_0|$ for the
    mitochondria example with functional given
    by~\eqref{eq:mitochondria:J} (to three significant
    digits). Computations performed with $T = 35$, $\kappa_n = 0.5$
    and $N_x = 40$ and the ESDIRK4 scheme.}
  \label{tab:mito:taylor}
\end{table}

In a multi-stage scheme, a number of intermediate stage solutions are
computed during the solution process. To avoid excessive memory usage,
dolfin-adjoint does not store these stage solutions. Thus, in order to
compute the adjoint solution, the stage solutions are recomputed for
each time-step. As a consequence, the minimal ratio of adjoint runtime 
to forward runtime for multistage ODE solves using this strategy
is approximately 2. For a
linear PDE solve, the optimal ratio of adjoint run time to forward run
time, assuming all linear solves take the same amount of time, is $1$. For a nonlinear PDE
solve via a Newton iteration, the optimal ratio of adjoint run
time to forward run time, assuming all linear solves take the same
amount of time, is
$1/N$ where N is the number of average Newton iterations in the
forward solve. These optimal ratios for linear and nonlinear PDE
solves have observed for typical dolfin-adjoint
usage~\cite{FarrellEtAl2013}. For this example, two Newton iterations
were required on average to solve the PDEs to within a tolerance of
$10^{-10}$ both for $N_x = 40$ and $80$, and assembly dominated the
PDE solve runtime. Thus for this example with two multistage ODE
solves and a nonlinear PDE solve in each timestep, the optimal ratio
of adjoint runtime to forward runtime is expected to be in the range
$0.5-2$, and closer to the upper bound than the lower bound depending
on the distribution of computational cost between the PDEs and ODEs.

Experimentally observed timings for this example are listed in
Tables~\ref{tab:mito:runtimes}--\ref{tab:mito:runtimes:odes} for $N_x
= 40, 80, 160$. From Table~\ref{tab:mito:runtimes}, we observe that
the adjoint-to-forward runtime ratio for the total solve is in the
range $1.38-1.95$ and decreasing with increasing problem size. The
corresponding gradient-to-forward runtime is in the range $1.44 -
2.28$, again decreasing with increasing problem size as expected. From
Table~\ref{tab:mito:runtimes:odes}, we observe that adjoint-to-forward
runtime ratio for the ODE solves is in the range $2.35-2.41$ for this
set of problem sizes, also decreasing with increasing problem size. We
also note that the ODE solve runtime, both for the forward and adjoint
solves, appears to scale linearly with the problem size $M = N_x^2$ as
anticipated and as is optimal. The ODE solves account for $37-40\%$ of
the total forward runtime for this example.
\begin{table}
  \centering
  \begin{tabular}{c|c|cc|cc}
    \toprule
    $N_x$ & Forward $t_F$ (s) & Adjoint $t_A$ (s) & Ratio $\frac{t_A}{t_F}$ & Gradient $t_G$ (s) & Ratio $\frac{t_G}{t_F}$ \\
    \midrule
    $40$  & 9.30 & 18.1 & 1.95 & 21.2 & 2.28 \\
    $80$  & 34.0 & 52.2 & 1.53 & 56.1 & 1.65 \\
    $160$ & 135  & 186  & 1.38 & 193 & 1.44 \\
    \bottomrule
  \end{tabular}
  \vspace{1em}
  \caption{Run times to compute the forward solution $t_F$ and to
    compute the adjoint $t_A$ and the functional gradient $t_G$ of the
    functional defined by~\eqref{eq:mitochondria:J} for the
    mitochondria example with increasing spatial resolution for $T =
    35$ $\kappa_n = 0.5$ with the ESDIRK4 multistage scheme.}
  \label{tab:mito:runtimes}
\end{table}
\begin{table}
  \centering
  \begin{tabular}{c|ccc|c}
    \toprule
    $N_x$ & Forward ODEs (s) & (of total) & Adjoint ODEs (s) & Ratio \\
    & $t_{F, O}$ &  & $t_{A, O}$ & $\frac{t_{A, O}}{t_{F, O}}$ \\
    \midrule
    $40$  & 3.48 & (37 \%) & 8.40 & 2.41 \\
    $80$  & 13.7 & (40 \%) & 32.6 & 2.37 \\
    $160$ & 53.6 & (40 \%) & 126  & 2.35 \\
    \bottomrule
  \end{tabular}
  \vspace{1em}
  \caption{Run times for solving the ODE systems when computing the
    forward solution $t_{F, O}$ and when computing the adjoint $t_{A,
      O}$ of the functional defined by~\eqref{eq:mitochondria:J} for
    $T = 35$ $\kappa_n = 0.5$ with the ESDIRK4 multistage scheme.}
  \label{tab:mito:runtimes:odes}
\end{table}

We also conducted numerical experiments varying the multistage scheme
used in the ODE solves, including the 1-stage Backward Euler (BDF1),
the 2-stage Crank-Nicolson (CN2), the 3-stage ESDIRK3 scheme in
addition to the previously reported 4-stage ESDIRK4 scheme. The
results (using $N_x = 40$) are reported in
Table~\ref{tab:mito:runtimes:odeschemes}. We observe that the
adjoint-to-forward runtime ratio ranges from $2.04$ (BDF1) to $2.41$
(ESDIRK3, ESDIRK4) and the general trend is that this ratio increases
with the number of stages as expected, but stabilizes around
$2.3-2.4$.
\begin{table}
  \centering
  \begin{tabular}{l|ccc|c}
    \toprule
    Scheme & Forward & Forward ODEs (s) & Adjoint ODEs (s) & Ratio \\
    & $t_F$ & $t_{F, O}$ & $t_{A, O}$ & $\frac{t_{A, O}}{t_{F, O}}$ \\
    \midrule
    BDF1 & 6.56 & 0.76 & 1.55 & 2.04 \\
    CN2 & 6.67 & 0.94 & 2.14 & 2.28 \\
    ESDIRK3 & 8.30 & 2.54 & 6.02 & 2.37 \\
    ESDIRK4 & 9.30 & 3.48 & 8.40 & 2.41 \\
    \bottomrule
  \end{tabular}
  \vspace{1em}
  \caption{Run times computing the forward solution $t_F$, for solving
    the ODE systems when computing the forward solution $t_{F, O}$ and
    when computing the adjoint $t_{A, O}$ of the functional defined
    by~\eqref{eq:mitochondria:J} for $T = 35$ $\kappa_n = 0.5$, $N_x =
    40$ for some common multistage schemes.}
  \label{tab:mito:runtimes:odeschemes}
\end{table}

\subsection{Application: cardiac electrophysiology (2D)}
\label{sec:cardiac:2D}

In this example, we consider the bidomain
equations~\eqref{eq:bidomain} over a two-dimensional rectangular
domain $\Omega = [0, 50] \times [0, 50]$ (mm) with coordinates ($x_0,
x_1$). We will consider a set of different cardiac cell models of
increasing complexity: a reparametrized FitzHugh-Nagumo (FHN) model
with 1 ODE state variable~\cite{FitzHugh1961}, the Beeler-Reuter (BR)
model with 7 ODE state variables~\cite{BeelerReuter1977}, the ten
Tusscher and Panfilov (TTP) epicardial cell model with 18 ODE state
variables~\cite{tenTusscherPanfilov2006b}, and the Grandi et al (GPB)
cell model with 38 ODE state variables~\cite{grandi2010novel}. The
description of each of these cell models is available via the CellML
repository, and our implementation of the ionic current $I_{\rm ion}$
and $F$ arising in~\eqref{eq:bidomain} was automatically generated
from the CellML models for the BR, TTP and GPB models. We refer to the
supplementary code~\citep{TheCode} for the precise description of the
models including our choice of FitzHugh-Nagumo parameters.

Our parameter setup is otherwise as follows. We let $\chi = 140$
(mm$^{-1}$), $C_m = 0.01$ ($\mu$F/mm$^2$), and let $M_i =
\mathrm{diag} (g_{\rm if}, g_{\rm if})$ and $M_e = \mathrm{diag}
(g_{\rm ef}, g_{\rm es})$ where $g_{\rm ef} = 0.625/(\chi C_m)$,
$g_{\rm es} = 0.236/(\chi C_m)$ and $g_{\rm if} = 0.174/(\chi
C_m)$. We let $I_s = 0$ and set the initial condition
\begin{equation}
  v_0 = 10 \left ( \frac{x_0}{50} \right ) ^2 + 10
\end{equation}
for $v$. The other state variables are initialized to the default
initial conditions given by the respective CellML models.

We use a second-order Strang splitting scheme (i.e.~$\theta = 0.5$
for~\eqref{eq:bidomain:odes}--\eqref{eq:bidomain:pdes}), a
Crank-Nicolson discretization in time and continuous piecewise linear
finite elements in space for the PDE system, and different Rush-Larsen
type discretizations (RL1, GRL1, RL2, GRL2) for the time
discretization of the ODE systems\footnote{We observed spurious wave
  propagation results when using a first-order splitting scheme and
  implicit Euler for with the BR model.}. The coupled linear systems
that arise at each PDE step~\eqref{eq:bidomain:pdes} were solved with
a block-preconditioned GMRES scheme with a relative tolerance of
$10^{-10}$. The complete details of the solver are detailed in the
supplementary code~\citep{TheCode}.

We take $\kappa_n = 0.1$ (ms), and consider a uniform mesh of the
computational domain with $N_x \times N_x \times 2$ triangles. The
coarsest mesh considered used $N_x = 160$. For this case,
simulations converged using any of the cell models and schemes, and
gave qualitatively correct results, except the Grandi cell model
discretized by the RL1 scheme for which the numerical solution scheme
failed to converge.

\begin{figure}
  \centering
  \begin{subfigure}[b]{0.32\textwidth}
    \includegraphics[width=\textwidth]{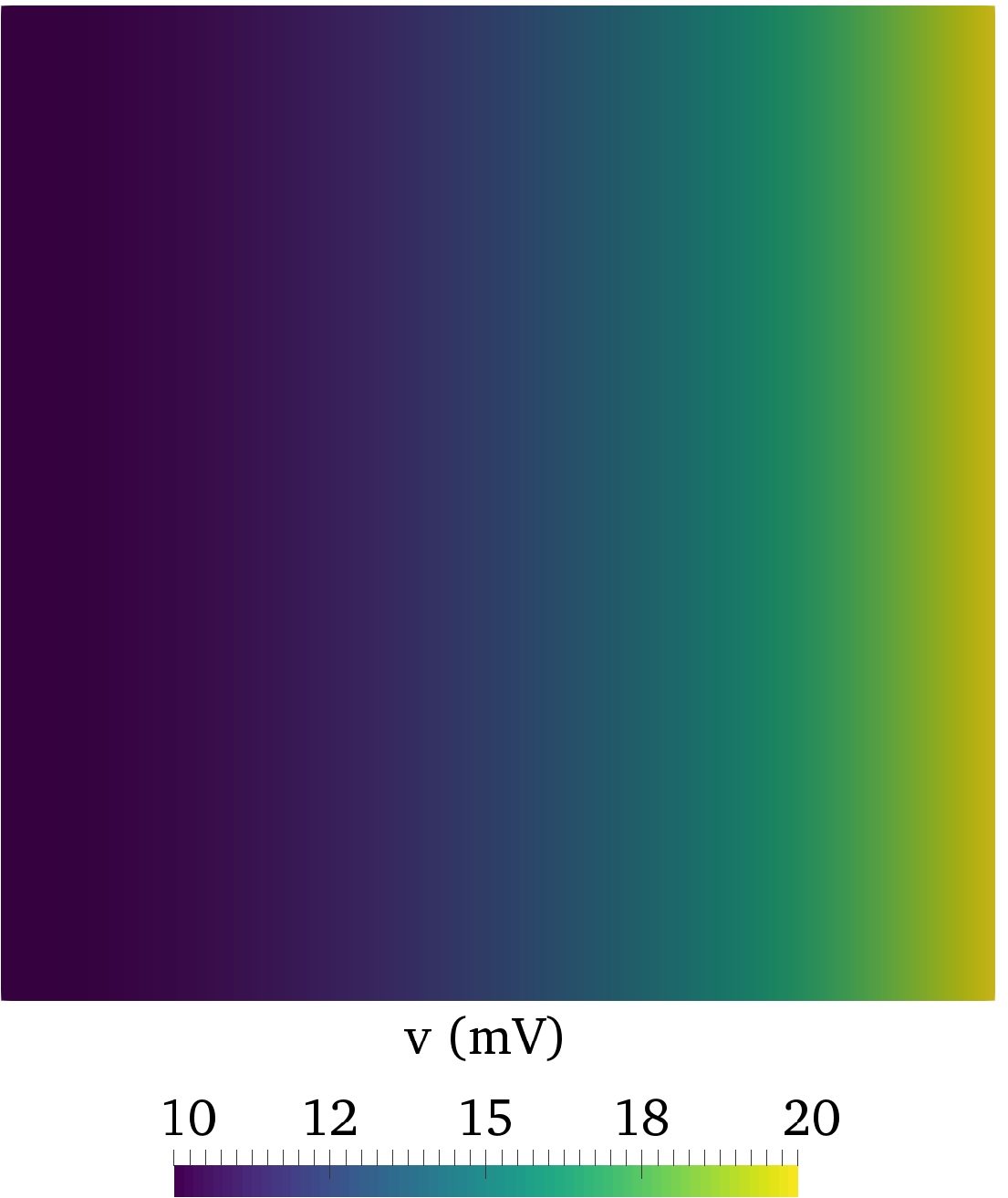}
    \caption{t=0 ms}
  \end{subfigure}
  \begin{subfigure}[b]{0.32\textwidth}
    \includegraphics[width=\textwidth]{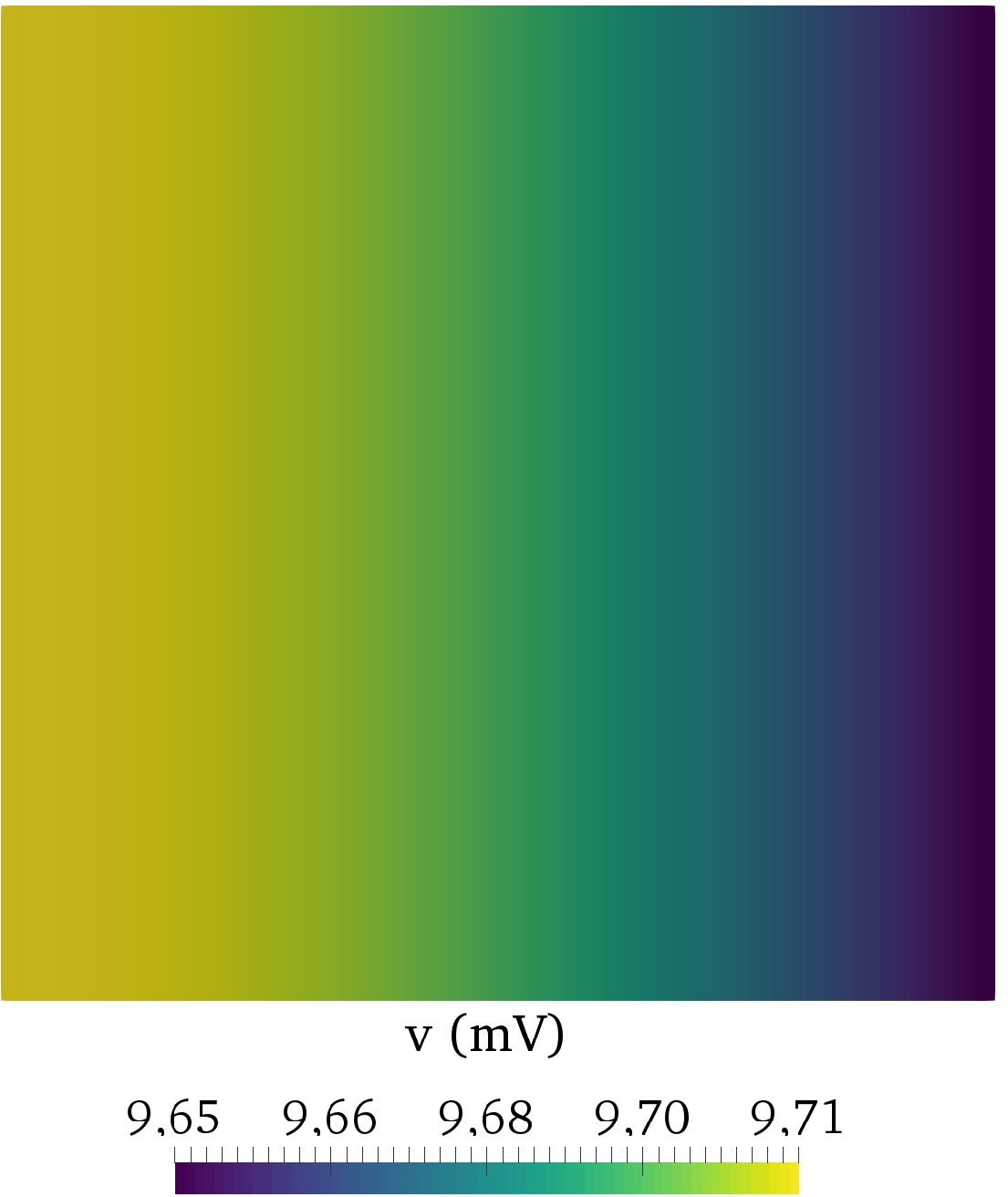}
    \caption{t=100 ms}
  \end{subfigure}
  \begin{subfigure}[b]{0.316\textwidth}
    \includegraphics[width=\textwidth]{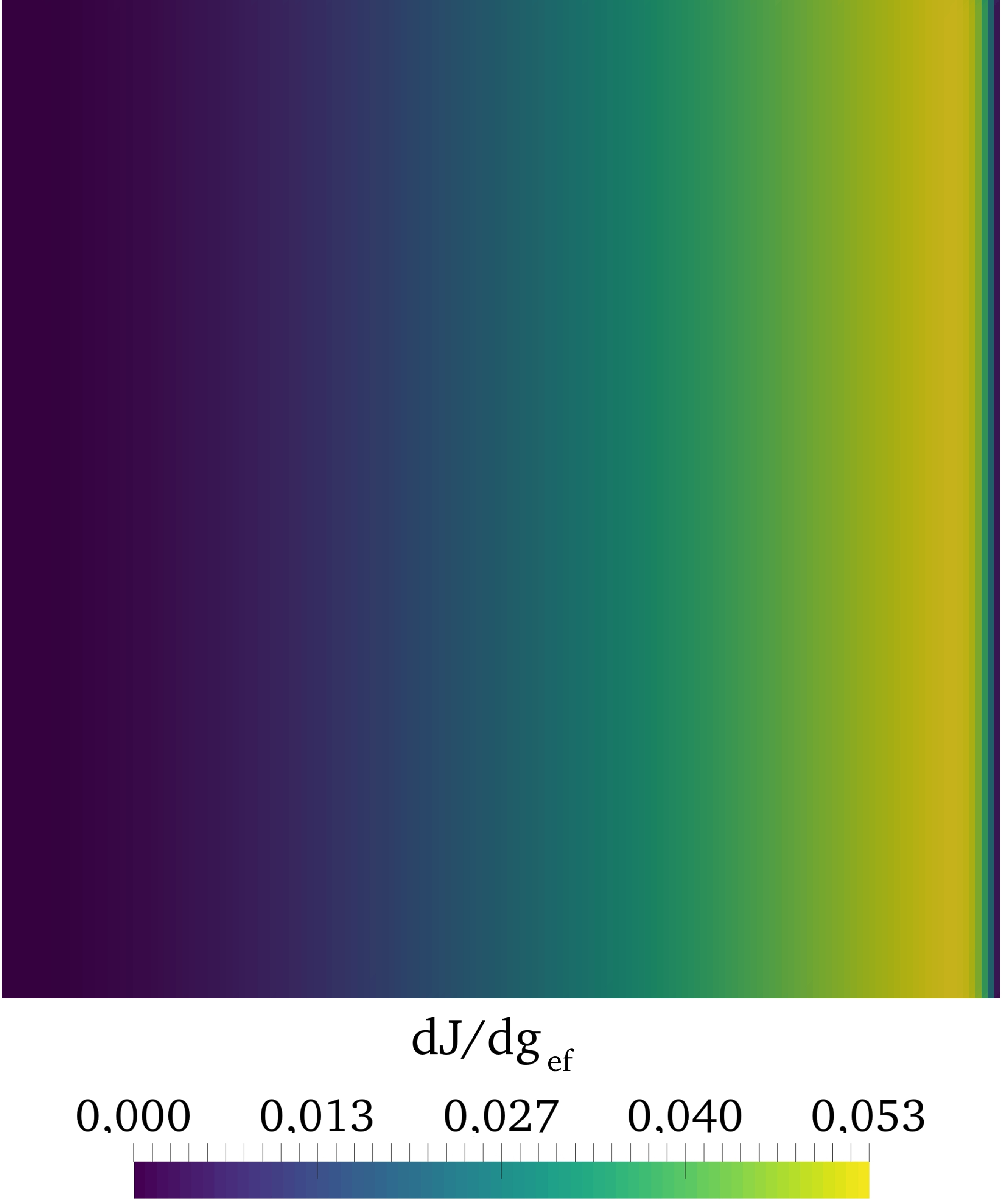}
    \caption{dJ/dg$_{ef}$}
  \end{subfigure}\\  
  \begin{subfigure}[b]{0.6\textwidth}
    \includegraphics[width=\textwidth]{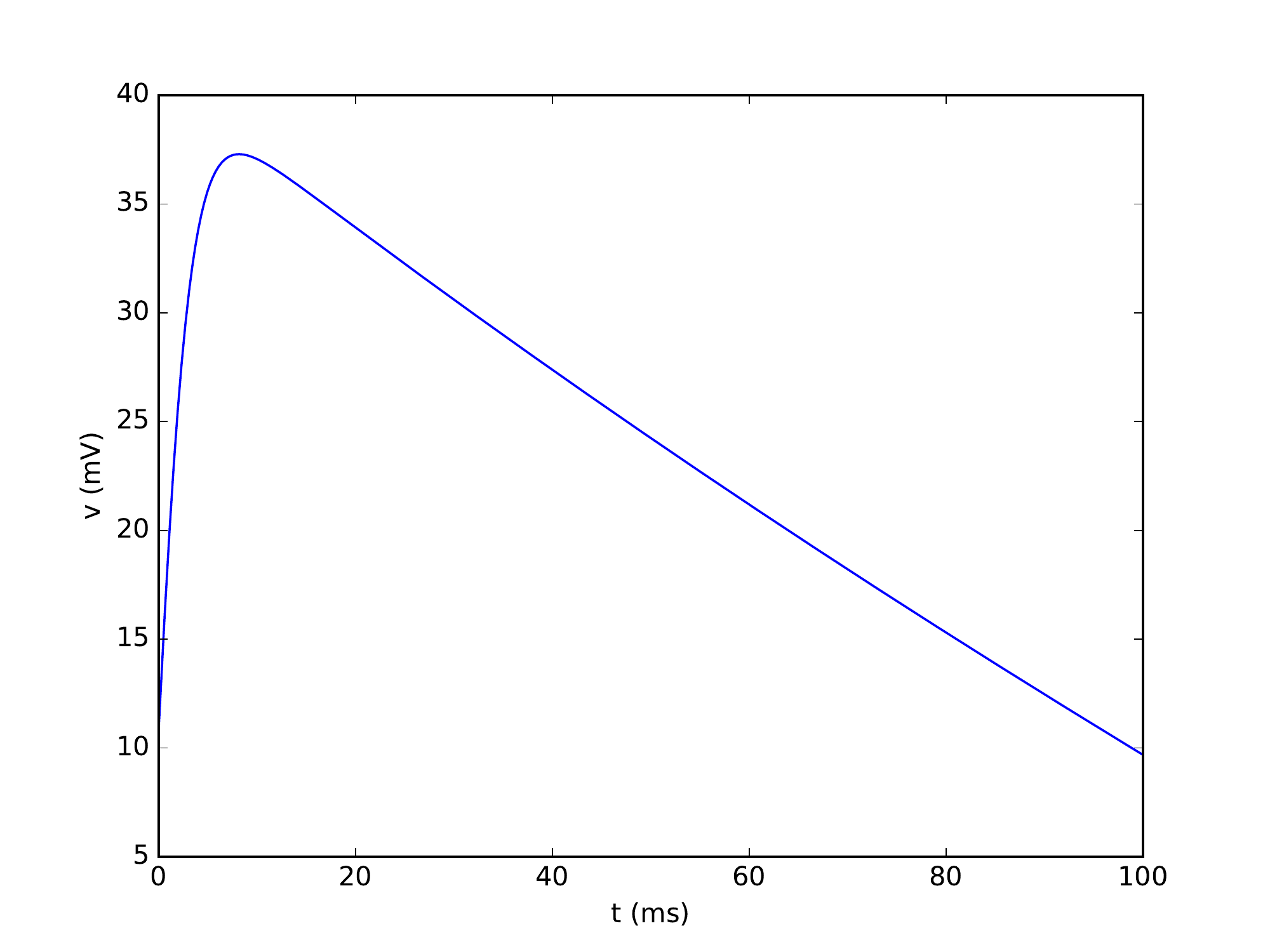}
    \caption{}
  \end{subfigure}  
  \caption{(a) Transmembrane potential $v$ at $t = 0$ ms and (b) at $t=100$ ms using the FHN model, computed with $\kappa_n=0.1$ and $N_x=160$. 
    (c) $L^2$-gradient of the functional with respect to the extracellular fiber conductivity $g_{ef}$. (d) Timeplot of transmembrane potential $v$ at $x = y = 25$ mm.}
  \label{fig:cardiac:setup}
\end{figure}

We are interested in computing the gradient with respect to the
extracellular fiber conductivity $g_{\rm ef}$ of the following
objective functional
\begin{equation}
  J(v) = \sum_{i=1}^5 \int_{\Omega} \left ( v(t_i)^2 + s(t_i)^2 \right ) \dx
  \label{eq:cardiac:J}
\end{equation}
for a equidistributed set of time points $t_i$. The model initial condition and solution for the
transmembrane potential at $T = 100$, together with the $L^2$-gradient
of the objective functional with respect to the fiber conductivity are
illustrated in Figure~\ref{fig:cardiac:setup}.

The following numerical experiments were performed on the Abel
supercomputer at the University of Oslo. All runtime experiments were
repeated three times, of which the minimum timing values are reported
here.

\subsubsection{Verifying the correctness of the discrete gradient}
\label{sec:verifying_gradient}

To verify that the computed functional gradient is correct, we
performed a Taylor test in a random perturbation direction at given
spatial and temporal resolutions. The results with the Fitz-Hugh
Nagumo cell model and the GRL1 scheme are listed in
Table~\ref{tab:cardiac:taylor_test} and demonstrate the expected
orders of convergence without a computed gradient (order 1) and when
using the computed discrete gradient (order 2) in the Taylor
expansion.
\begin{table}
  \centering
  \begin{tabular}{l|cccc}
    \toprule
    $s$ & $R_0 (s)$ & order & $R_1 (s)$ & order \\
    \midrule
    0.005   & 15.7  & 0.96 & 0.541    & \\
    0.0025  & 7.97  & 0.98  & 0.145    & 1.90 \\
    0.00125 & 4.02  & 0.99 & 0.0382   & 1.92 \\
    0.000625 & 2.02  & 0.99 & 0.00983   & 1.96 \\
    \bottomrule
  \end{tabular}
  \vspace{1em}
  \caption{Taylor convergence results for the cardiac
    electrophysiology 2D test case. Taylor remainders $R_0 = |J(g_{ef}
    + s\delta)|$ and $R_1 = |J(g_{ef} + s\delta) - J(g_{ef}) - \nabla
    J(g_{ef}) s \delta|$ for decreasing step sizes $s$ and a random
    direction $\delta$ with functional given
    by~\eqref{eq:cardiac:J}. Computations performed with $T = 10.0$,
    $\kappa_n = 0.1$ and $N_x = 160$, the FHN cell model and the GRL1
    scheme. We observe that the remainders converge at first order
    without gradient information and at second order with gradient
    information, as expected.}
  \label{tab:cardiac:taylor_test}
\end{table}
We also performed similar experiments for the other cell models (BR,
TTP, and GPB) and other Rush-Larsen solution schemes (RL1, GRL2,
RL2). We used the same settings as in
Table~\ref{tab:cardiac:taylor_test}, except for the Grandi cell model
with any other Rush-Larsen scheme than GRL1, where the time step had to
be reduced to 0.01 in order for the forward solver to converge. We
obtained the expected convergence order for all combination of cell
models and schemes tested. These results indicate that the
automatically derived and computed adjoints and gradient are correct
for all the Rush-Larsen schemes implemented.

\subsubsection{Adjoint runtime performance}\label{sec:adjoint_performance}

As the bidomain PDEs~\eqref{eq:bidomain:pdes} are linear, a rough
theoretical estimate of the optimal PDE adjoint to PDE forward runtime
ratio is 1, assuming that the solution times in the adjoint and
forward PDE solves are equal and dominated by the finite element
assembly and linear solvers. Similarly, the ODEs are solved with
explicit Rush-Larsen schemes with a theoretically optimal
adjoint-to-forward ratio of 1 as well.

In Table~\ref{tab:adjoint_to_forward}, we list the total runtime, and
the runtime components corresponding to only the ODE solves, only the
PDE solves and an intermediate variable update (merge) step, for both
the forward and the adjoint solution computation for a test case with
$T = 10.0$, $\kappa_n = 0.1$, $N_x = 160$, the GRL1 scheme and the
four different cardiac cell models (FHN, BR, TTP, GPB).

\begin{table}[h]
  \centering
  \begin{tabular}{llccccc}
    \toprule  
    Model &  & Total & ODEs & PDEs & Merge & Other \\
    \midrule
      FHN & Forward & 38.8      & 3.20 & 34.1 & 0.466 & 3 \% \\      
          & Adjoint & 44.4      & 8.27 & 28.1 & 1.75 & 14 \% \\      
          & Ratio   & 1.15        & 2.58  & 0.82 & 3.76 &   \\        
      \midrule   
      BR & Forward & 48.3       & 10.9 & 34.5 & 0.979 & 4  \% \\    
          & Adjoint & 62.1      & 20.4 & 28.2 & 4.66  & 14 \% \\     
          & Ratio   & 1.29      & 1.88 & 0.81 & 4.76  &  \\    
      \midrule   
      TTP & Forward & 67.3      & 24.4 &  35.5  & 1.97 &  8   \% \\
          & Adjoint & 99.5      & 47.8 &  28.2  & 9.98 &  14  \% \\
          & Ratio   & 1.48      & 1.96 &  0.79  & 5.07 &         \\    
      \midrule   
      GPB & Forward & 91.8      & 40.8 & 34.7   & 3.63 & 14   \% \\
          & Adjoint & 152       & 81.7 & 28.2   & 19.6 & 15   \% \\
          & Ratio   & 1.66      & 2.00 & 0.81   & 5.40 &         \\    
    \bottomrule
  \end{tabular}
  \caption{Adjoint-to-forward runtime performance for a cardiac
    electrophysiology 2D test case. Breakdown of forward and adjoint
    runtimes and their ratio for different cell models. Each row shows
    total runtime (s), runtime for the ODE, PDE and merge steps (s),
    and percentage of the total runtime not accounted for by these
    three steps (Other). The computations were performed with $T =
    10$, $\kappa_n = 0.1$, $N = 160$, GRL1 as the ODE scheme and
    \eqref{eq:cardiac:J} as the adjoint functional.}
  \label{tab:adjoint_to_forward}
\end{table}

From Table~\ref{tab:adjoint_to_forward}, for the forward solves, we
observe that the total runtime increases with cell model
complexity. The forward PDE step runtime is essentially independent of
the cell model, while the forward ODE step runtime increases with the
cell model complexity, as expected: the ODE runtime corresponds to
approximately $10 \%$ of the PDE runtime for the simplest
FitzHugh-Nagumo model, and approximately $114 \%$ of the PDE runtime
for the most complex Grandi model. The runtime of the forward merge
step is insignificant in comparison to the ODE and PDE steps, but
increases with cell model complexity.

We observe that also the adjoint PDE runtimes are comparable for all
cell models. The resulting adjoint-to-forward PDE step ratios are in
the range $0.79-0.82$. The decrease in compute time for the adjoint
PDE solves compared to the forward PDE solves are attributable to the Krylov
solvers: fewer Krylov solver iterations were necessary for convergence
for the adjoint solves compared to the forward solves.
For the ODE solves, we observe that the adjoint-to-forward ratio for
the ODE solves is in the range $1.88-2.58$ for the different cell
models. Additionally, we observe that the runtime of the adjoint merge step is
significantly larger than the corresponding time for the forward merge
step and with an increasing adjoint-to-forward ratio with increasing
cell model complexity $3.76-5.40$. We suggest that the reason for the
increased and increasing adjoint runtime is that the adjoint merge
step involves additional assembly of variational forms of complexity
comparable to that of the cell model.

Overall, we observe that the adjoint-to-forward total runtime ratios
are in the range $1.15-1.66$ for the different cell models considered,
with increasing adjoint-to-forward ratio with increasing cell model
complexity. Since the PDE runtime dominates the total runtime for the
simpler cell models, the total adjoint-to-forward ratio is closer to
$1$ for these models. On the other hand, for the more complicated cell
models, the ODE runtimes dominate, and so we observe
adjoint-to-forward ratios closer to $2$ for these.

In the last column of Table~\ref{tab:adjoint_to_forward}, we observe
that the combined ODE and PDE solve time contributes to between 97\%
(for FitzHugh-Nagumo) and 86\% (for Grandi) of the total forward
runtime for the cell models considered. The remaining forward runtime
cost is primarily caused by the initialization routines such as
loading the computational mesh, and creating the required function
spaces. For the corresponding adjoint runtimes, the combined ODE and
PDE solve times contribute to between 86\% (for FitzHugh-Nagumo) and
85\% (for Grandi) of the total run time, and the percentage
contribution decreases with the complexity of the cell model.  The
remaining runtime cost can primarily be attributed to recording the
forward states during the adjoint solve (which is included in the
reported total run time but not in the ODE or PDE runtimes here). We
also timed the computation of the discrete gradient and noted that the
cost of additionally computing the gradient was negligible (less than
0.1\% of the adjoint solve runtime).

\subsubsection{Forward and adjoint parallel scalability}
\label{sec:scaling}

In this section, we discuss the weak and strong parallel scalability
of the forward and adjoint ODE solvers and the scalability of the
point integral solvers in particular applied to this 2D
electrophysiology test case. We focus on the scalability of the point
integral solver particular as this is the primary new feature
described in this paper. Previous studies have discussed the parallel
scalability of FEniCS in general~\citep{RichardsonWells2015,
  AlnaesEtAl2015, FarrellEtAl2013}.

We ran a series of weak scaling experiments on Abel, hosted by the
Norwegian national computing infrastructure (NOTUR), for the
previously considered series of cell models (FHN, BR, TTP, GPB) on an
increasing number of CPUs (1, 4, 16, 64) and with an increasing mesh
resolution.  We timed the total runtime of the point integral solves,
both for the forward solves and the adjoint solves. For each cell
model and resolution, we ran three experiments and extracted the
minimal run times. The results are listed in
Table~\ref{tab:weak_scaling_odes}. For both the forward and adjoint
run times, we computed the parallel efficiency (PE) as the ratio of
the 1-CPU runtime to the $n$-CPU run time for $n = 4, 16, 64$, also
listed in Table~\ref{tab:weak_scaling_odes}. The optimal PE would be
1. The last column of the Table lists the computed the
adjoint-to-forward ratio for these solves.

For the forward runtimes, we observe that the parallel efficiencies
range from 0.66 to 0.95 for all cell models and resolutions, with
efficiencies higher than 0.8 for all but the least complex cell
model. In general, we observe that the parallel efficiency increases
with increasing cell model complexity and that it decreases moderately
with the number of CPUs. Some of the loss of parallel efficiency (from
1) is attributable to only a near-perfect distribution of mesh
vertices across CPUs; we observed approximately 10\% variation between
CPUs in the number of local vertices.

For the adjoint runtimes, we observe parallel efficiencies ranging
from 0.75 (FHN on 64 cores) to 0.95, again with efficiencies higher
than 0.8 for all but the least complex cell model. In general, we note
the same trends as for the forward run times: increasing parallel
efficiency with increasing cell model complexity and moderately
decreasing parallel efficiency with increasing number of CPUs. We
emphasize that we thus achieve the \emph{same parallel efficiency} for
the adjoint point integral solves as for the forward point integral solves.

Comparing the forward and adjoint runtimes, we see that the
adjoint-to-forward ratio for the point integral solvers range from
1.19 to 1.56. Comparing these numbers with the adjoint-to-forward
ratios of the overall ODE solves as discussed in the previous, we
observe that the adjoint-to-forward point integral solver performance
is better. The adjoint-to-forward ratio is moderately increasing with
cell model complexity, but is stable with respect to the number of
CPUs. This latter point again illustrates that the adjoint point
integral solvers demonstrate the same parallel efficiency as the
forward point integral solves.
\begin{table}
  \centering
  \begin{tabular}{lccccccc}
    \toprule
     & \#CPUs & $N_x$ & Forward PIS (s) & PE & Adjoint PIS (s) & PE & Ratio\\
    \midrule
FHN & 1  & 160         & 0.185    &       &   0.251  &      & 1.36 \\
    & 4  & 320         & 0.253    & 0.73  &   0.311  & 0.81 & 1.23 \\
    & 16 & 640         & 0.263    & 0.70  &   0.330  & 0.76 & 1.25 \\
    & 64 & 1280        & 0.281    & 0.66  &   0.333  & 0.75 & 1.19 \\
\midrule    
BR & 1    & 160         & 0.740    &       & 0.925   &      & 1.25 \\
   & 4    & 320         & 0.795    & 0.93  & 0.998   & 0.93 & 1.26 \\
   & 16   & 640         & 0.868    & 0.85  & 1.076   & 0.86 & 1.24 \\
   & 64   & 1280        & 0.899    & 0.82  & 1.143   & 0.81 & 1.27 \\
\midrule    
TTP & 1     & 160       & 1.657   &      & 2.516     &      & 1.52 \\
    & 4     & 320       & 1.736   & 0.95 & 2.640     & 0.95 & 1.52 \\
    & 16    & 640       & 1.886   & 0.88 & 2.939     & 0.86 & 1.56 \\
    & 64    & 1280      & 1.991   & 0.83 & 3.072     & 0.82 & 1.54 \\
\midrule 
GPB & 1     & 160       & 2.618   &       & 3.895    &      & 1.49 \\
    & 4     & 320       & 2.771   & 0.94  & 4.131    & 0.94 & 1.49 \\
    & 16    & 640       & 3.059   & 0.86  & 4.553    & 0.86 & 1.49 \\
    & 64    & 1280      & 3.125   & 0.84  & 4.632    & 0.84 & 1.48 \\
    \bottomrule
  \end{tabular}
  \caption{Weak scaling of the point integral solver for a cardiac
    electrophysiology 2D test case. Forward and adjoint point integral
    solver runtimes for increasing $N_x$ and simultaneously increasing
    number of cores and for different cell models. Each row shows
    number of CPUs, the mesh resolution, the total point integral
    solver runtime (s) for the forward solves, the parallel efficiency
    (PE) for the point integral forward solves, the total point
    integral solver runtime (s) for the adjoint solves, the parallel
    efficiency for the point integral adjoint solves, and the
    adjoint-to-forward runtime ratio. The computations were performed
    with $T = 1.0$, $\kappa_n = 0.1$ and GRL1 as the ODE scheme and
    \eqref{eq:cardiac:J} as the adjoint functional.}
  \label{tab:weak_scaling_odes}
\end{table}

We also ran a series of strong scaling experiments on Abel for the
previously considered series of cell models (FHN, BR, TTP, GPB) on an
increasing number of CPUs (1, 4, 16, 64) and with a fixed mesh
resolution. We timed the total runtime of the point integral solver
steps for both the forward and adjoint solves. Again, for each cell
model and resolution, we ran three experiments and extracted the
minimal run times. The results are listed in
Table~\ref{tab:strong_scaling}. For both the forward and adjoint run
times, we computed the (strong scaling) parallel efficiency (PE) as
the ratio of the 1-CPU runtime to $n$-times the $n$-CPU runtime for $n
= 4, 16, 64$, also listed in Table~\ref{tab:strong_scaling}. The
optimal (strong scaling) PE would be 1. The last column of the Table
lists the computed the adjoint-to-forward ratio for these solves.

We observe that the parallel efficiencies are comparable to the
results from the weak scaling test case with efficiencies in the range
from 0.60 to 0.94. We note that the adjoint strong parallel efficiency
is comparable to the forward strong parallel efficiency. The
adjoint-to-forward ratios for this strong scaling test are also in the
same range as for the weak scaling test.

\begin{table}
  \centering
  \begin{tabular}{lccccccc}
    \toprule
      Model  & No. cores   & Forward PIS & PE          & Adjoint PIS & PE   & Ratio\\
      \midrule 
      FHN    &     $1$     &  12.874     &             &  17.039     &      & 1.32 \\ 
             &     $4$     &  3.503      & 0.92        &   4.826     & 0.88 & 1.38 \\
             &     $16$    &  1.159      & 0.69        &   1.373     & 0.78 & 1.18 \\
             &     $64$    &  0.281      & 0.72        &   0.333     & 0.80 & 1.19 \\
    \midrule                         
      BR     &      $1$    &  48.286     &             &  59.673     &      & 1.24 \\
             &      $4$    &  13.760     & 0.88        &  16.455     & 0.91 & 1.20 \\
             &      $16$   &   3.525     & 0.86        &   4.307     & 0.87 & 1.22 \\
             &      $64$   &   0.899     & 0.84        &   1.143     & 0.82 & 1.27 \\
    \midrule                      
      TTP    &     $1$     & 105.451     &             &  161.443    &      & 1.53 \\
             &     $4$     &  28.003     & 0.94        &  44.613     & 0.90 & 1.59 \\
             &     $16$    &  8.323      & 0.79        &  11.451     & 0.88 & 1.38 \\
             &     $64$    &   1.991     & 0.83        &   3.072     & 0.82 & 1.54 \\
    \midrule                         
      GPB    &     $1$     &  166.582    &             &  249.980    &      & 1.50 \\
             &     $4$     &   45.811    & 0.91        &   68.084    & 0.92 & 1.49 \\
             &     $16$    &   17.426    & 0.60        &   25.787    & 0.61 & 1.48 \\
             &     $64$    &    3.125    & 0.83        &   4.632     & 0.84 & 1.48 \\
    \bottomrule
  \end{tabular}
  \caption{Strong scaling of the point integral solver for a cardiac
    electrophysiology 2D test case. Forward and adjoint point integral
    solver runtimes for a fixed (fine) mesh resolution $N_x = 1280$
    with an increasing number of CPUs, for different cell models (FHN,
    BR, TTP, GPB). Each row show the number of CPUs, the total point
    integral solver runtime for the forward solves, the corresponding
    parallel efficiency (PE), the total point integral solver runtime
    for the adjoint solves, the corresponding parallel efficiency (PE)
    and the adjoint-to-forward runtime ratio. The computations were
    performed with $T = 1.0$, $\kappa_n = 0.1$ and GRL1 as the ODE scheme
    and \eqref{eq:cardiac:J} as the adjoint functional.}
  \label{tab:strong_scaling}
\end{table}

\subsection{Application: biventricular cardiac electrophysiology (3D)}
\label{sec:cardiac:3D}

In this example, we aim to compute the sensitivity of the squared
$L^2$-norm of the transmembrane potential with respect to its initial
condition, in order to illustrate the features discussed in this paper
for a more complex test case. This example also aims to illustrate how
different representations may seamlessly be used for the control
variables; this is useful for instance when computing sensitivities
with respect to spatially constant, regionally defined or highly
resolved control variables.

Inspired by~\citep{ArevaloEtAl2016}, we consider the monodomain
variation of~\eqref{eq:bidomain} over a mesh of a three-dimensional
bi-ventricular domain $\Omega$ with the ten Tusscher and Panfilov
epicardial cell model~\citep{tenTusscherPanfilov2006b} as detailed in
Section~\ref{sec:cardiac:2D}. For the monodomain equations, we replace
the bidomain equations~\eqref{eq:bidomain:1}--\eqref{eq:bidomain:2}
by: find $v$ such that
\begin{equation}
  \label{eq:monodomain}
  v_t - \Div M \Grad v = - I_{\rm ion}(v, s)  \quad \text{ in } \Omega.
\end{equation}
with homogeneous Neumann boundary conditions and with $v(t = 0) =
v_0$. As before, we let $\chi = 140$ (mm$^{-1}$), $C_m = 0.01$
($\mu$F/mm$^2$), and for simplicity (ignoring realistic,
spatially-varying fiber directions), let $M = \mathrm{diag} (g_{\rm f},
g_{\rm s}, g_{\rm n})$ where:
\begin{equation}
  g_f = 0.255, \quad g_n = 0.0775, \quad g_s = 0.0775.
\end{equation}

\begin{figure}
  \centering
  \includegraphics[width=0.48\textwidth]{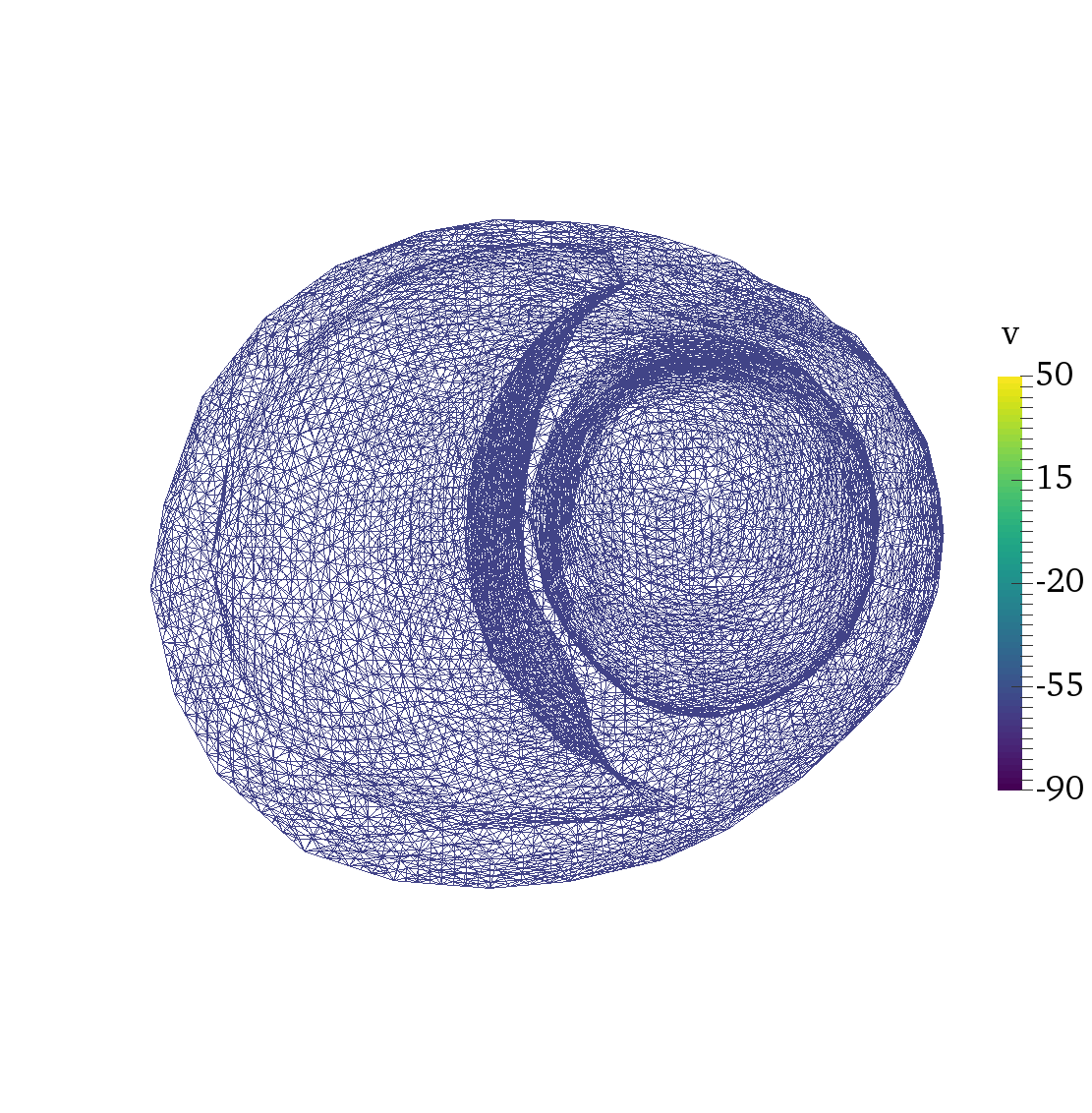}
  \includegraphics[width=0.48\textwidth]{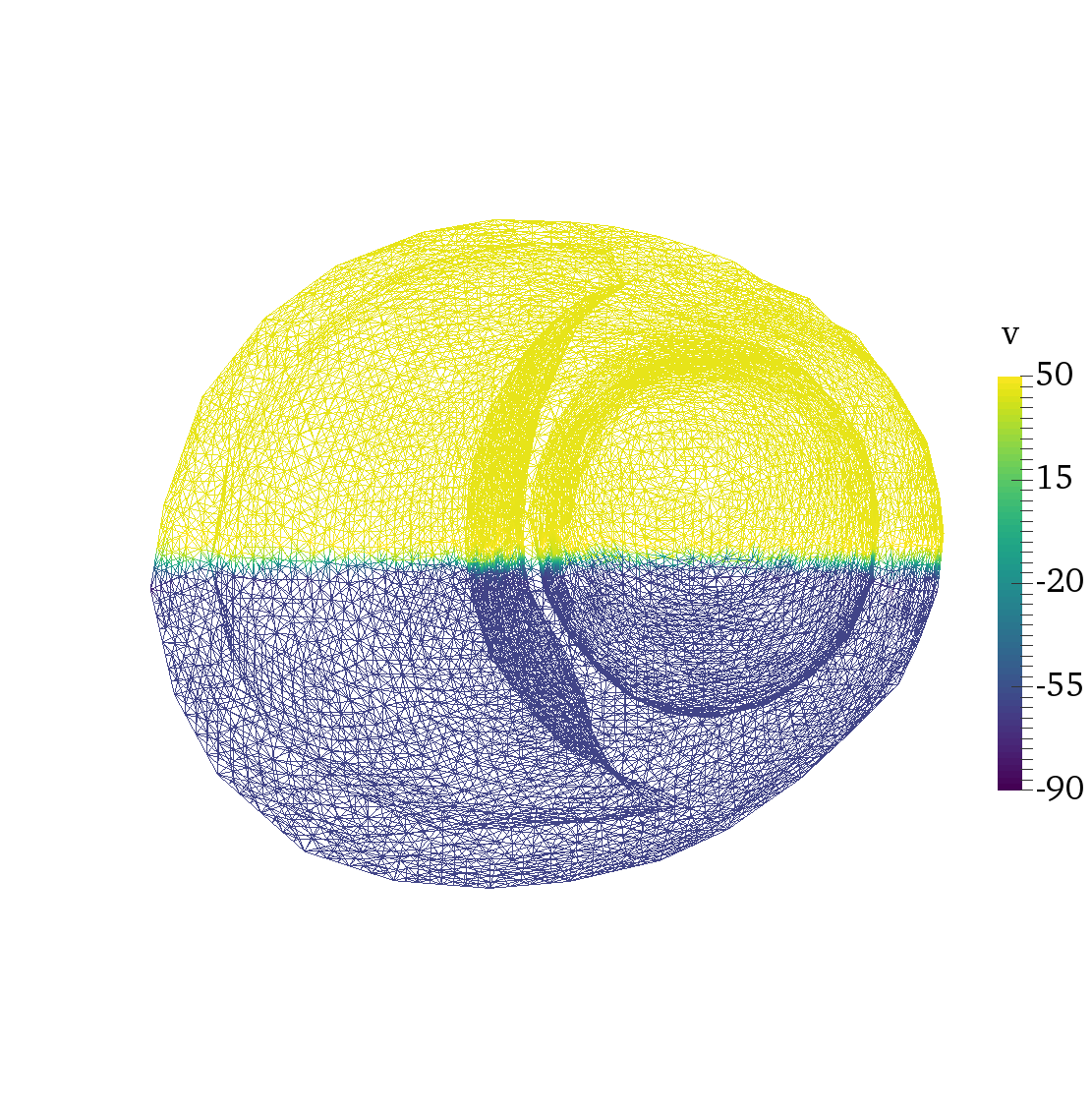}
  \caption{Transmembrane potentials: wireframe view (along positive
    x-axis) of the computed transmembrane potential $v(t)$ at $t = 0$
    (left) and $t = 6.0$ ms (right).}
  \label{fig:cardiac:3D:solutions}
\end{figure}

To solve the coupled equations~\eqref{eq:monodomain}
and~\eqref{eq:bidomain:odes}, we use a second-order Strang splitting scheme
($\theta = 0.5$),
a Crank-Nicolson discretization in time and continuous piecewise
linear finite elements in space for the resulting PDEs, and the
first-order generalized Rush-Larsen (GRL1) scheme for the resulting
ODEs. We take $\kappa_n = 0.05$ (ms). The moderately coarse mesh has
mesh cell diameters in the range $[0.61, 2.70]$ (mm) and bounding box
$[0.103, 40.2] \times [-15.8, 36.8] \times [-30, 13.9]$ (mm$^3$), with
45 625 vertices and 236 816 cells, and is illustrated in
Figure~\ref{fig:cardiac:3D:solutions}.

To quantify the transmembrane potential exceeding a zero threshold ,
we consider the following objective functional:
\begin{equation}
  J(v) = \int_{\Omega_v} v(T)^2 \dx,
  \label{eq:cardiac3D:J}
\end{equation}
at a fixed time $T = 6.0$ (ms) and where $\Omega_v = \{ x \in \Omega
\, | \, v(x) \geq 0 \}$. We are
interested in computing the gradient of $J$ with respect to the initial
condition $u_0$.

We consider two different representations of the initial
condition to illustrate computing the gradient with respect to low and
high resolution fields. For both cases, we let
\begin{equation}
  v_0(x) = \begin{cases}
    -61.1 & x_2 > z_{\rm mid}, \\
    -61.3 & x_2 \leq z_{\rm mid}, \\
  \end{cases}
  \label{eq:cardiac3D:v0}
\end{equation}
with $z_{\rm mid} = -8.07$, but consider (i) $v_0$ as represented by
two constants $(c_0, c_1)$ i.e.~the two-dimensional space $\R^2$ and
(ii) $v_0$ represented by a continuous piecewise linear defined
relative to the mesh i.e. as a $n_{\rm verts}$-dimensional space $V$
where $n_{\rm verts}$ denotes the number of mesh vertices.

The initial condition and computed solution at $t = 6.0$ (ms) are
illustrated in Figure~\ref{fig:cardiac:3D:solutions}. (With reference
to the chosen perspective in this Figure, we denote the subdomain
where $x_2 > z_{\rm mid}$ as the \emph{top} part of the domain and the
subdomain $x_2 \leq z_{\rm mid}$ as the \emph{bottom} part.)

\begin{figure}
  \centering
  \includegraphics[width=0.45\textwidth]{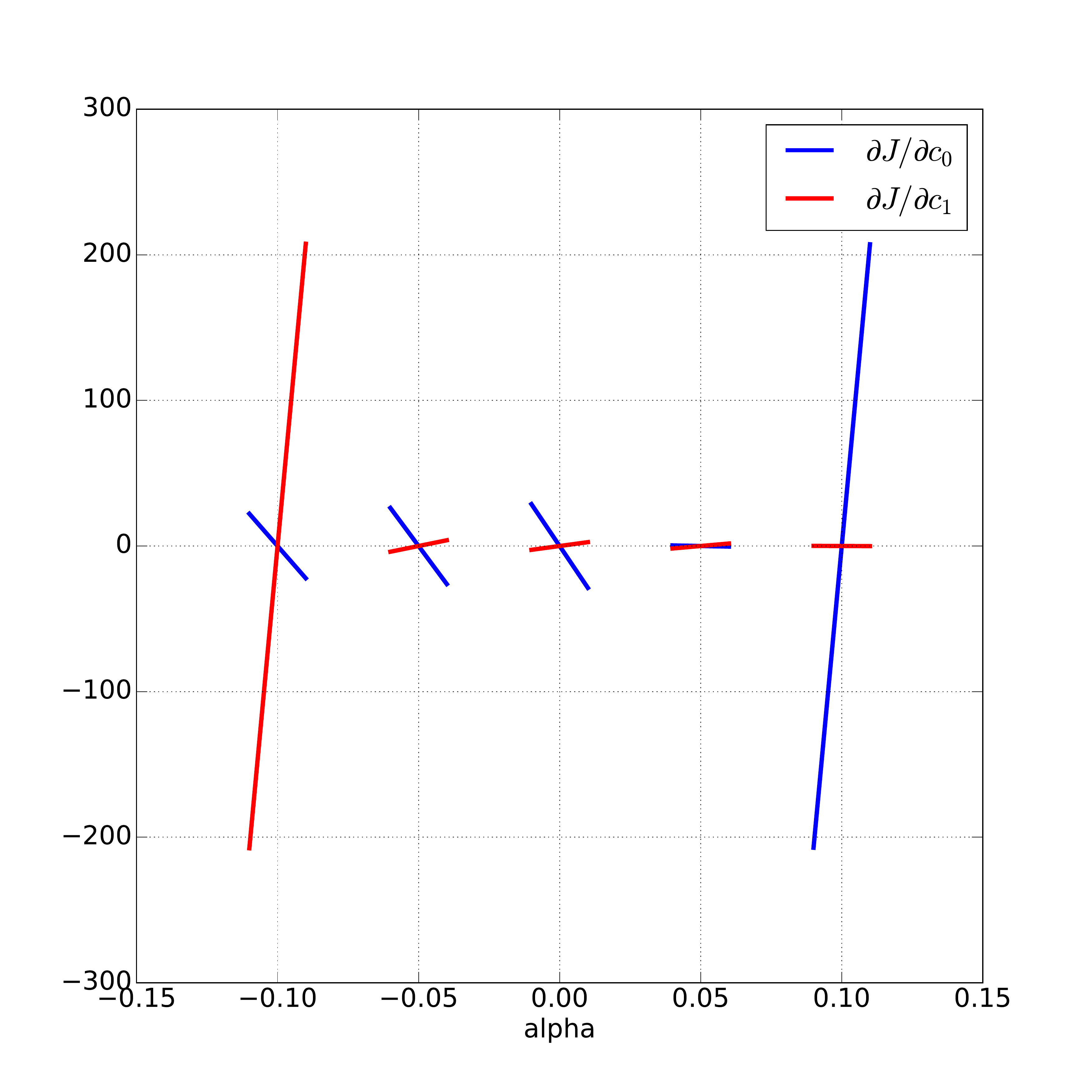}
  \hspace{-2em}
  \includegraphics[width=0.55\textwidth]{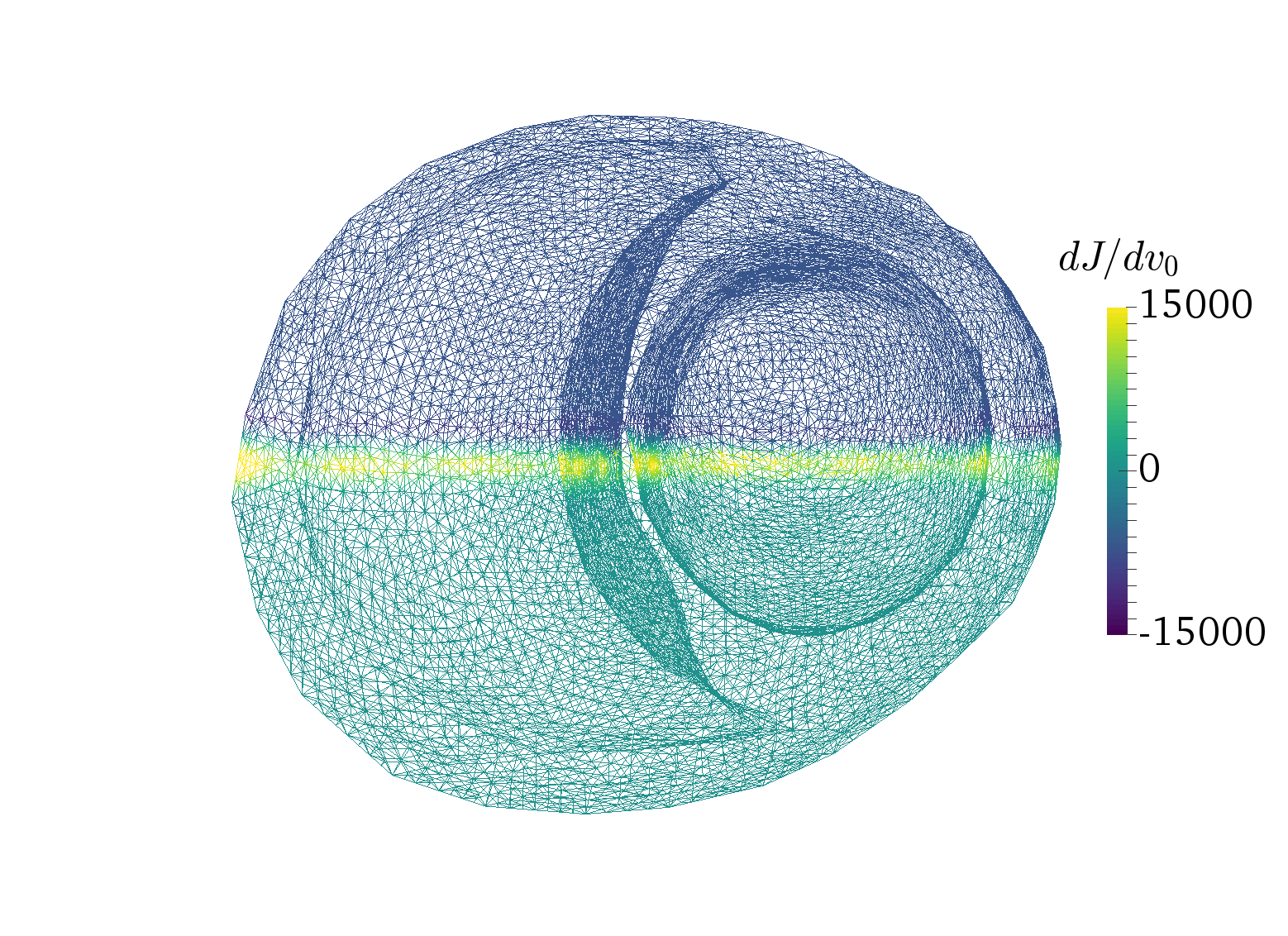}
  \caption{Sensitivities computed with respect to low and high
    resolution representations of an initial transmembrane
    potential. Left: sensitivities $\tfrac{\partial J}{\partial v_0} =
    ( \tfrac{\partial J}{\partial c_0}, \tfrac{\partial J}{\partial
      c_1} )$ of the objective functional $J$ defined
    by~\eqref{eq:cardiac3D:J} with respect to a low resolution initial
    condition $v_0 = (c_0, c_1) \in \R^2$ for different initial
    conditions $v_0 = v_0^{\alpha} = (-61.1 + \alpha, -63.3 + \alpha)$
    for $\alpha = \{ -0.1, -0.05, 0.0, 0.05, 0.1 \}$. For each point
    $\alpha$ and $i = 0, 1$, each plotted line segment (spanning from
    $\alpha - 0.01$ to $\alpha + 0.01$) has slope $\tfrac{\partial
      J}{\partial c_i}$.  Right: Sensitivity $\tfrac{\partial
      J}{\partial v_0}$ of the objective functional $J$ defined
    by~\eqref{eq:cardiac3D:J} with respect to high resolution initial
    condition $v_0 \in V$.}
  \label{fig:cardiac:3D:gradient}
\end{figure}

For the case where the initial condition is represented at a low
resolution by $v_0 = (c_0, c_1) \in \R^2$, we considered a series of
experiments starting with different initial conditions $v_0^{\alpha} =
(c_0 + \alpha, c_1 + \alpha) \in \R^2$ for $\alpha \in \{ -0.1, -0.05,
0.0, 0.05, 0.1 \}$. The resulting gradients (in $\R^2$) with respect
to $v_0 \in \R^2$ are illustrated in
Figure~\ref{fig:cardiac:3D:gradient} (left). We observe that as we
increase $\alpha$, the gradient with respect to the top initial
condition $c_0$ goes from large and positive to small and negative,
while the gradient with respect to the bottom initial condition $c_1$
varies from moderately negative to close-to-zero and rapidly to large
and positive. The rapid changes in the gradients illustrate the
nonlinear nature of the problem at hand.

The sensitivity of the objective functional $J$ defined
by~\eqref{eq:cardiac3D:J} with respect to the initial condition $v_0
\in V$ is illustrated in Figure~\ref{fig:cardiac:3D:gradient}
(right). We observe that the gradient in the upper part of the domain
is large and negative, that the gradient in a boundary zone is large
and positive and the gradient in the lower part of the domain is close
to zero. This indicates that infinitesimal increases/decreases in the
upper part of the domain will lead to a large decrease/increase in the
functional value, infinitesimal increases/decreases in the boundary
zone domain will lead to a large increase/decrease in the functional
value, and that infinitesimal increases/decreases in the lower part of
the domain will have relatively small effect on the functional value.

\section{Concluding remarks}
\label{sec:conclusion}

We have presented a high-level framework that allows for the
specification of coupled PDE-ODE systems and their efficient forward
and adjoint solution with operator splitting schemes. We have
illustrated the features of the framework with a series of examples
originating from cell modelling and computational cardiac
electrophysiology. Our numerical results indicate adjoint runtime
performance near optimal (measured in terms of adjoint-to-forward
runtimes), and parallel efficiency indices of $84-94\%$ for the more
realistic cardiac cell models.

The main advantages of the framework are the speed of developing
solvers for new systems, and the ability to automatically derive their
adjoints from the high-level specification. This allows for flexible
exploration of models when the relevant governing equations are
uncertain, and for the identification of unknown parameters via the
solution of inverse problems as demonstrated
e.g.~in~\citep{KallhovdEtAl2017}. Both of these capabilities are of
significant importance in numerous areas of scientific computing, and
in particular in computational medicine, physiology and
biology. Future work will focus on further improving the parallel
scalability of the solvers and its application to problems in
personalized medicine.

\bibliographystyle{siamplain}
\bibliography{references}

\end{document}